\documentclass[final]{siamltex}

\usepackage{amssymb}
\usepackage{amsmath}
\usepackage{graphicx}
\usepackage{latexsym}
\usepackage{verbatim}
\usepackage{float}
\usepackage{graphics}
\usepackage{url}
\usepackage{subfig}
\usepackage{helvet}

\usepackage[english]{babel}
\usepackage[T1]{fontenc}
\usepackage[ansinew]{inputenc}
\usepackage{lmodern}	
\usepackage{amsmath}
\usepackage{amsfonts}
\usepackage{tikz}

\newtheorem{Thm}{Theorem}

\definecolor{DarkGreen}{rgb}{0,.55,0}

\title{Numerical Methods for the Computation of the Confluent and Gauss Hypergeometric Functions\thanks{This work was supported by the Numerical Algorithms Group (NAG) and the Engineering and Physical Sciences Research Council (EPSRC).}}

\author{John W. Pearson \thanks{School of Mathematics, Statistics and Actuarial Science, University of Kent, Cornwallis Building (East), Canterbury, Kent, CT2 7NF, UK ({\tt j.w.pearson@kent.ac.uk})} \and Sheehan Olver \thanks{School of Mathematics and Statistics, The University of Sydney, Australia ({\tt Sheehan.Olver@sydney.edu.au})} \and Mason A. Porter \thanks{Oxford Centre for Industrial and Applied Mathematics, Mathematical Institute, University of Oxford, Radcliffe Observatory Quarter, Woodstock Road, Oxford, OX2 6GG, UK ({\tt porterm@maths.ox.ac.uk})}}

\begin{document}
\maketitle


\vspace{.5 in}


\begin{abstract}
The two most commonly used hypergeometric functions are the confluent hypergeometric function and the Gauss hypergeometric function.  We review the available techniques for accurate, fast, and reliable computation of these two hypergeometric functions in different parameter and variable regimes. The methods that we investigate include Taylor and asymptotic series computations, Gauss-Jacobi quadrature, numerical solution of differential equations, recurrence relations, and others. 
We discuss the results of numerical experiments used to determine the best methods, in practice, for each parameter and variable regime considered. 
We provide `roadmaps' with our recommendation for which methods should be used in each situation. 
\end{abstract}






\pagestyle{myheadings}
\thispagestyle{plain}
\markboth{J. W. PEARSON, S. OLVER, AND M. A. PORTER}{COMPUTATION OF HYPERGEOMETRIC FUNCTIONS}


\section{Introduction} \label{sec:intro}

The aim of this review paper is to summarize methods for computing the two most commonly used hypergeometric functions: the confluent hypergeometric function $\mathbf{M}(a;b,z)$ and the Gauss hypergeometric function $\mathbf{F}(a,b;c;z)$. (We also consider the associated functions $_{1}F_{1}(a;b;z)$ and $_{2}F_{1}(a,b;c;z)$.) We overview methods that have been developed for computing $\mathbf{M}$ and $\mathbf{F}$, and we discuss how to choose appropriate methods for different parameter and variable regimes. We thereby obtain reliable and fast computation for a large range of the parameters ($a$ and $b$ for $\mathbf{M}$; and $a$, $b$, and $c$ for $\mathbf{F}$) and the variable $z$.  Because accurate error bounds are seldomly available, we test a large variety of approaches, which we require to be stable, accurate, fast, and robust within the parameter and variable regions for which they have been selected. This is especially important when working with finite-precision arithmetic.

The computation of confluent hypergeometric functions and Gauss hypergeometric functions is important in a wide variety of applications \cite{Seaborn}.  For instance, these functions arise in areas such as photon scattering from atoms \cite{Gavrila}, networks \cite{PastorSatorrasVespignani,TorrieriValenti}, Coulomb wave functions \cite{BellScott,EUTENAEH, NobleThompson}, binary stars \cite{PPS}, mathematical finance \cite{BoylePotapchik},  non-Newtonian fluids \cite{ZhaoYang} and more.

Except for specific situations, computing hypergeometric functions is difficult in practice. A plethora of methods exist for computing each hypergeometric function; these include Taylor series, asymptotic expansions, continued fractions, recurrence relationships, hyperasymptotic expansions, and more.  However, each method is only reliable and efficient for limited parameter and variable regimes.  Consequently, one must be prepared to use different members of this suite of possibilities in different regimes of parameter and variable values.

We have also developed {\scshape Matlab} code for computing the functions $\mathbf{M}$ and $\mathbf{F}$ using a range of methods; this code is available at \cite{PearsonCodeLink}.


\section{Background on Hypergeometric Functions} \label{sec:background}

The confluent hypergeometric function is defined by \cite{DLMF}
\begin{align} \label{1F1} 
	\mathbf M(a;b;z) &= \sum_{j=0}^{\infty}\frac{(a)_{j}}{\Gamma(b+j)}~\frac{z^{j}}{j!}\,,
\end{align}
where the \textbf{Pochhammer symbol} $(\mu)_{j}$ is
\begin{align*}
	 (\mu)_{0}=1\,,\quad(\mu)_{j}=\mu(\mu+1) \times \cdots \times (\mu+j-1)\,,\quad j=1,2, \ldots ~\,.
\end{align*}
The function $\mathbf M$ is entire (i.e., it is analytic throughout the complex plane $\mathbb C$) in  the parameters $a$ and $b$ and the variable $z$. Therefore, the sum \eqref{1F1} always converges.  When $b$ is a non-positive integer, we define 
\begin{equation*}
	_{1}F_{1}(a;b;z)  =  \Gamma(b) \mathbf M(a;b;z) = \sum_{j=0}^{\infty}\frac{(a)_{j}}{(b)_{j}}~\frac{z^{j}}{j!}\,,
\end{equation*}
which is also commonly denoted by $M(a;b;z)$ and is itself often referred to as the confluent hypergeometric function.  Because $_{1}F_{1}(a;b;z)$ is not defined if $b$ is equal to a non-positive integer and there are numerical issues in its computation if $b$ is close to a non-positive integer, it is preferable to compute $\mathbf M$ when possible.

The Gauss hypergeometric function is defined within the unit disk $|z|<1$ by \cite{DLMF}
\begin{align} \label{2F1} 
	\mathbf F(a,b;c;z) &= \sum_{j=0}^{\infty}\frac{(a)_{j}(b)_{j}}{\Gamma(c+j)}~\frac{z^{j}}{j!}
\end{align}
and defined outside of the unit disk by analytic continuation, with a principal branch cut along $[1,\infty)$.  On $[1,\infty)$, it is defined to be continuous from below; in other words, $\mathbf F(a,b;c;z) = \lim_{\epsilon \rightarrow 0} \mathbf F(a,b;c;z-|{\epsilon}| i)$.  When $c$ is a non-positive integer, we define
\begin{align} 
	_{2}F_{1}(a,b;c;z) & = \Gamma(c)  \mathbf F(a,b;c;z) \qquad\left(= \sum_{j=0}^{\infty}\frac{(a)_{j}(b)_{j}}{(c)_{j}}~\frac{z^{j}}{j!}  \quad\hbox{for}~|z| < 1\right).
\end{align}
The function $_{2}F_{1}(a,b;c;z)$ is commonly denoted by $F(a,b;c;z)$ and is also frequently referred to as the Gauss hypergeometric function.  Analogous to the case for the confluent hypergeometric function, it is better numerically to compute $\mathbf F$  than $_{2}F_{1}$.



The confluent hypergeometric function $_{1}F_{1}(a;b;z)$ satisfies the ordinary differential equation \cite{Slater1}
\begin{align} \label{1f1de} 
	z\frac{{\rm d}^{2}f}{{\rm d}z^{2}}+(b-z)\frac{{\rm d}f}{{\rm d}z} - af=0\,,
\end{align}
while the Gauss hypergeometric function $_{2}F_{1}(a,b;c;z)$ satisfies \cite{Slater2}
\begin{align} \label{2f1de} 
	z(1-z)\frac{{\rm d}^{2}f}{{\rm d}z^{2}}+[c-(a+b+1)z]\frac{{\rm d}f}{{\rm d}z}-abf=0\,,
\end{align}
provided that none of $c$, $c-a-b$ or $a-b$ are equal to an integer. We discuss the numerical solution of these equations in Sections \ref{sec:1F1_othermethods} and \ref{sec:2F1_othermethods}, respectively.

The confluent and Gauss hypergeometric functions are examples of \textbf{special functions}, whose theory and computation are prominent in the scientific literature \cite{AbramowitzStegun,AAR,Hochstadt,LozierOlver,Luke4,Luke5,Moshier,Olver1,Temme3,Temme5,ZhangJin}. A large number of common mathematical functions (including $e^{z}$, $(1-z)^{-a}$, and $\sin^{-1}z$ \cite{AbramowitzStegun,Luke4}) and a large number of special functions (such as the modified Bessel function, the (lower) incomplete gamma and error functions, Chebyshev and Legendre polynomials, and other families of orthogonal polynomials \cite{AbramowitzStegun,AncaraniGasaneo,Mathar,Roach}) can be expressed in terms of confluent or Gauss hypergeometric functions.



\begin{figure}
  \centering
\subfloat{\label{fig:hyperegs1}\includegraphics[width=0.5\textwidth]{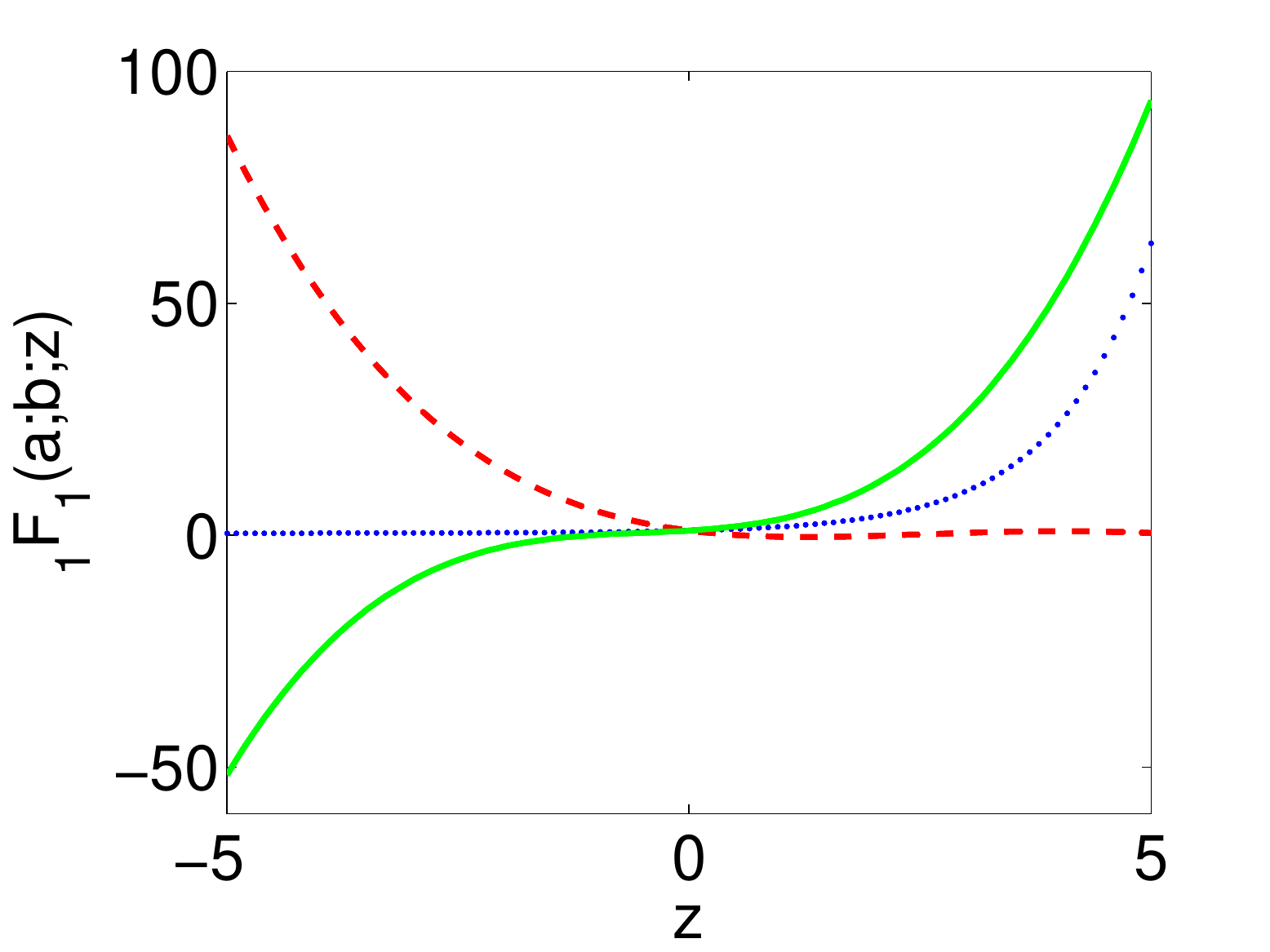}}
\subfloat{\label{fig:hyperegs2}\includegraphics[width=0.5\textwidth]{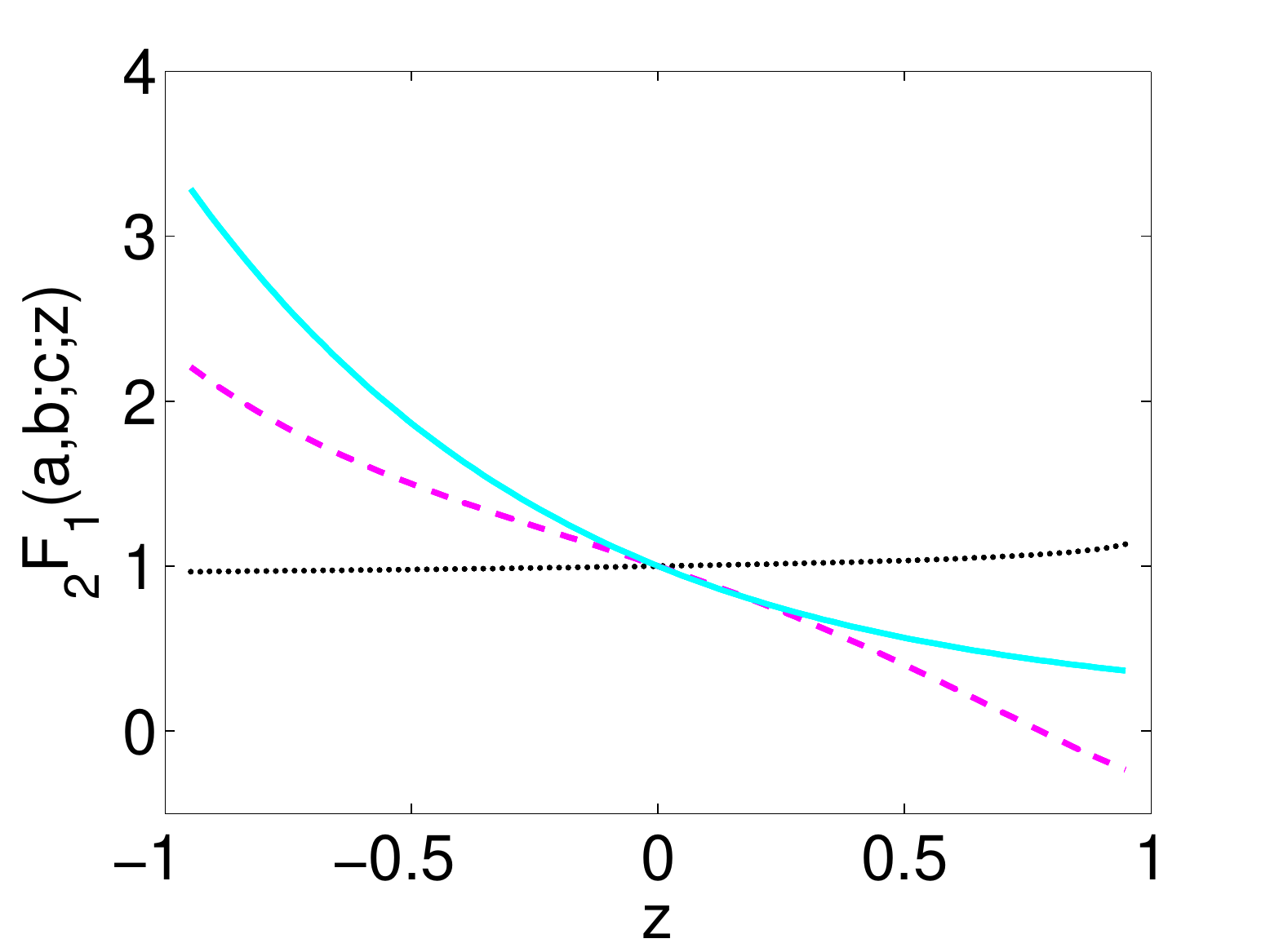}}
\caption{(Color online) (Left) Plots of $_{1}F_{1}(a;b;z)$, generated using {\scshape Matlab} \cite{Matlab}, for real $z\in[-5,5]$ with parameter values $(a,b)=(0.1,0.2)$ in dark blue (dotted), $(a,b) = (-3.8,1.5)$ in red (dashed), and $(a,b) = (-3,-2.5)$ in green (solid).  (Right) Plots of $_{2}F_{1}$ for real $z\in[-1,1]$ with parameter values $(a,b,c)=(0.1,0.2,0.4)$ in black (dotted), $(a,b,c) = (-3.6,-0.7,-2.5)$ in purple (dashed), and $(a,b,c) = (-5,1.5,6.2)$ in sky blue (solid).
}
\label{fig:hyperegs}
\end{figure}

\begin{figure}
  \centering
\subfloat{\label{fig:hyperegs3}\includegraphics[width=0.5\textwidth]{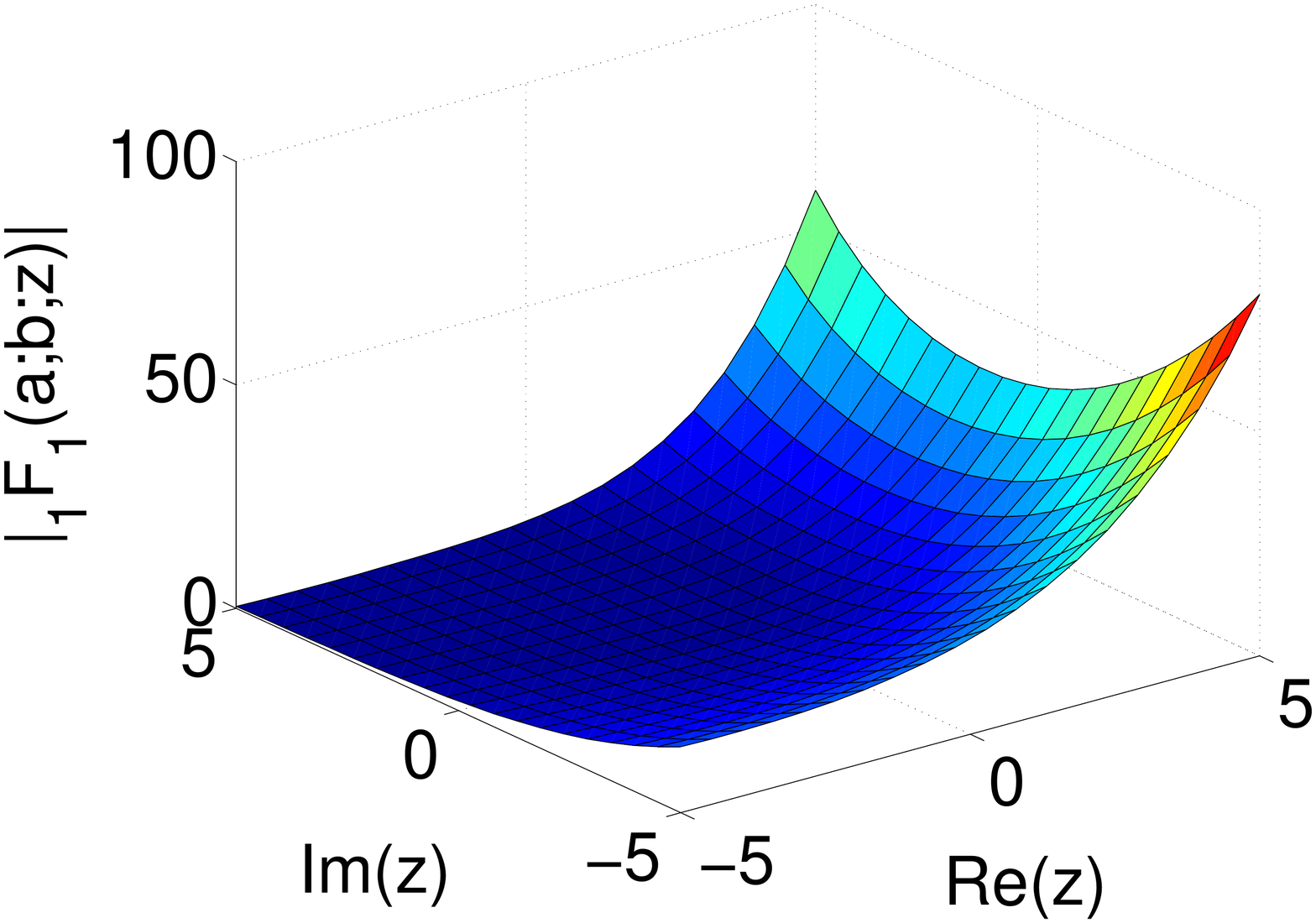}}
\subfloat{\label{fig:hyperegs4}\includegraphics[width=0.5\textwidth]{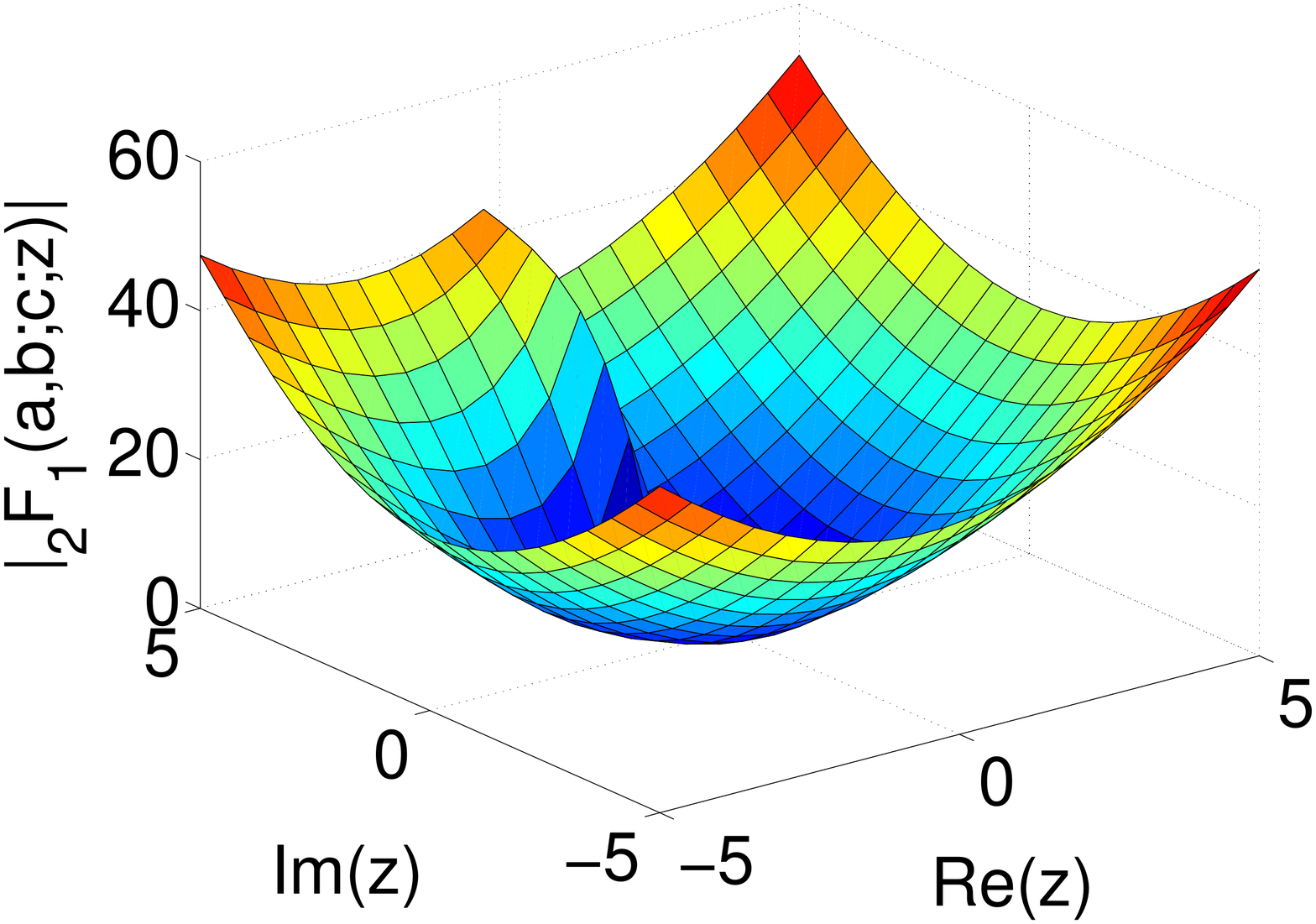}}
\caption{(Color online) (Left) Plot of $|_{1}F_{1}(a;b;z)|$ for $\big({\rm Re}(z),{\rm Im}(z)\big)\in[-5,5]^2$ with parameter values $(a,b)=(-3-2i,2.5)$. (Right) Plot of $|_{2}F_{1}|$ for $\big({\rm Re}(z),{\rm Im}(z)\big)\in[-5,5]^2$ with parameter values $(a,b,c)=(1+2i,-1.5,0.5-i)$.
}
\label{fig:hyperegscomplex}
\end{figure}

In Figs.~\ref{fig:hyperegs} and \ref{fig:hyperegscomplex}, we show plots of hypergeometric functions. These provide examples of the possible shapes of such functions in the real and complex planes, though of course there are a wide range of possible representations of them.

As we indicated above, the computation of hypergeometric functions arises in a wide variety of applications.  This underscores the importance of determining which computational techniques will be accurate and efficient for different parameter and variable regimes. It is also important to supply information on how to test the reliability of a routine, provide test cases that a routine might have difficulty computing (see Appendix \ref{sec:appendixa}), and indicate how to evaluate other special functions required for the computation of hypergeometric functions (see Appendix \ref{sec:appendixc}). As illustrated in Ref.~\cite{Pearson}, the inbuilt routines to compute hypergeometric functions in commercial software packages such as {\scshape Matlab} R2013a \cite{Matlab} and Mathematica 8 \cite{Mathematica} are not without limitations, so one should not rely on them in many parameter and variable regimes.  We have implemented our routines for {\scshape Matlab} R2013a, which uses double-precision arithmetic, and we 
have made it freely available online at \cite{PearsonCodeLink}. Of course, it is easier to obtain accurate results when computing these functions in higher-precision arithmetic, as is done when the inbuilt hypergeometric function commands in {\scshape Matlab} and {\scshape Mathematica} are called.  Our codes may be modified to run in higher precision.


\section{Computation of the Confluent Hypergeometric Function $\mathbf{M}$} \label{sec:1F1}

In this section, we discuss the methods that perform the best in practice for accurately and efficiently evaluating the confluent hypergeometric function. The range of methods that we discuss includes series computation methods (Sections \ref{sec:1F1_taylor}, 
\ref{sec:1F1_buchholz}, and \ref{sec:1F1_asymptotic}), quadrature methods (Section \ref{sec:1F1_quadrature}), recurrence relations (Section \ref{sec:1F1_recurrences}), and other methods (Section \ref{sec:1F1_othermethods}). In Section \ref{sec:1F1_summary}, we summarize the  computational strategies and provide recommendations for the most effective approaches for a variety of parameter regimes.


\subsection{Properties of $\mathbf{M}$} \label{sec:1F1_properties}


As stated in Section \ref{sec:background}, the \emph{confluent hypergeometric function} $\mathbf{M}(a;b;z)$ is defined by the series \eqref{1F1} for any $a\in\mathbb{C}$, $b\in\mathbb{C}$ and $z\in\mathbb{C}$. We note that $\mathbf M(a;b;0)={1 \over \Gamma(b)}$.  Additionally, if $a=n\in\mathbb{Z}^{-}\cup\{0\}$, then this series is given by a polynomial in $z$ of degree $-n$.


A second, linearly-independent solution of equation (\ref{1f1de}) is denoted $U(a;b;z)$ and is defined by the property that $U(a;b;z)\sim z^{-a}$ as $\left|z\right|\rightarrow\infty$ for $\left|\arg z\right|\leq\frac{3}{2}\pi-\delta$,
with $0<\delta\ll1$ \cite{Luke4}. The function $U(a;b;z)$ satisfies \cite{AbramowitzStegun}
\begin{align*}
	U(a;b;z) &= \frac{\pi}{\sin(\pi{}b)}\left(\frac{1}{\Gamma(1+a-b)} \mathbf M(a;b;z)-\frac{z^{1-b}}{\Gamma(a)}\mathbf M(1+a-b;2-b;z)\right)\,.
\end{align*}
The function $U(a;b;z)$ is defined with a principal branch cut along $(-\infty,0]$.  When $m\in\mathbb{Z}$, one obtains \cite{DLMF}
\begin{align*}
 	U(a;b;ze^{2\pi im})=\frac{2\pi ie^{-\pi ibm}\sin(\pi bm)}{\Gamma(1+a-b)\sin(\pi b)}~\mathbf{M}(a;b;z)+e^{-2\pi ibm}U(a;b;z)\,.
\end{align*}

 $\mathbf M(a;b;z)$ and $U(a;b;z)$ are also related for $b\notin\mathbb{Z}$ by \cite{Temme6}
\begin{align}
\label{bfMU}	\mathbf M(a;b;z) &= \frac{e^{\mp a\pi i}}{\Gamma(b-a)}U(a;b;z)+\frac{e^{\pm(b-a)\pi i}}{\Gamma(a)}e^{z}U(b-a;b;e^{\pm\pi i}z)\,.
\end{align}
Therefore, methods for computing $U(a;b;z)$ may also be useful for computing $\mathbf M$, and vice-versa, providing care is taken with numerical issues such as cancellation. One can also exploit methods for computing $_{1}F_{1}$ and $_{2}F_{1}$ in conjunction with the above expressions.

There are many useful applications of the confluent hypergeometric function. For example, it has been used to find exact solutions of the wave equation \cite{Hochstadt}, obtain insights on the scattering of charged particles \cite{BellScott}, analyze gene frequency data \cite{VGDB}, and investigate Asian options in finance \cite{BoylePotapchik}.


\subsection{Taylor Series} \label{sec:1F1_taylor}

The simplest method for computing the confluent hypergeometric function is to truncate the Taylor series
\begin{align} \label{1F1defnew} 
	_{1}F_{1}(a;b;z) \approx S_N =  \sum_{j=0}^{N}\underbrace{\frac{(a)_{j}}{(b)_{j}}~\frac{1}{j!}z^{j}}_{A_{j}}\,,
\end{align}
and then use the relation $_{1}F_{1}(a;b;z)=\Gamma(b)\mathbf{M}(a;b;z)$.

A \textbf{stopping criterion} must  be specified. It is common to stop computing terms when $\frac{\left|A_{N+1}\right|}{\left|S_{N}\right|}<\verb+tol+$ for some tolerance $\verb+tol+$ and some $N$, and then return $S_{N}$ (our approximation of $_{1}F_{1}(a;b;z)$) as the solution \cite{Muller}. Alternatively, one can terminate computation of the series using the more stringent condition that two successive terms are small. We use this latter approach in our numerical experiments.

By taking $\verb+tol+=\verb+eps+$, we obtain an accuracy of 15 decimal places (if round-off error does not play a role).  Note that ${\tt eps}\approx2.2\times10^{-16}$ denotes ``machine epsilon'' in {\scshape Matlab}. We choose this value for \url{tol} for all of the series computations in this paper.

Even the simple task of summing Taylor series can be accomplished in different ways, and this can have a significant impact on stability and efficiency.  We consider two methods for computing $S_N$.

\vspace{0.75em}
\underline{\textbf{Method (a):}} ~Compute
\begin{align*}
\ A_{0}{}&=1\,,{}\quad S_{0}{}=A_{0}\,, \\
\ A_{j+1}{}&=A_{j}~\frac{a+j}{b+j}~\frac{z}{j+1}\,,{}\quad S_{j+1}{}=S_{j}+A_{j+1}\,,\quad{}j=0,1,2,~\ldots\,.
\end{align*}

\vspace{0.75em}
\underline{\textbf{Method (b):}} ~The following three-term recurrence relation can be used to obtain approximations of $\mathbf{M}(a;b;z)$ recursively in terms of previous approximations \cite{Muller}:
\begin{align*}
\	S_{-1}&{}=S_{0}=1\,,\quad\quad\quad~{}S_{1}{}=1+\frac{a}{b}z\,, \\
\ r_{j}&{}=\frac{a+j-1}{j(b+j-1)}\,,\quad{}S_{j}=S_{j-1}+(S_{j-1}-S_{j-2})r_{j}z,{}&j{}=2,3,\ldots~\,.
\end{align*}
\vspace{0.75em}
%


As one can see in Table \ref{appbtable1} in Appendix \ref{sec:appendixb}, methods (a) and (b) generate similar results and require the same number of terms for most computations. However, as also illustrated by Table \ref{appbtable1}, method (b) is in general more effective when carrying out computations involving small parameters.

\begin{figure}[h]
\centering
\includegraphics[width=8.5cm]{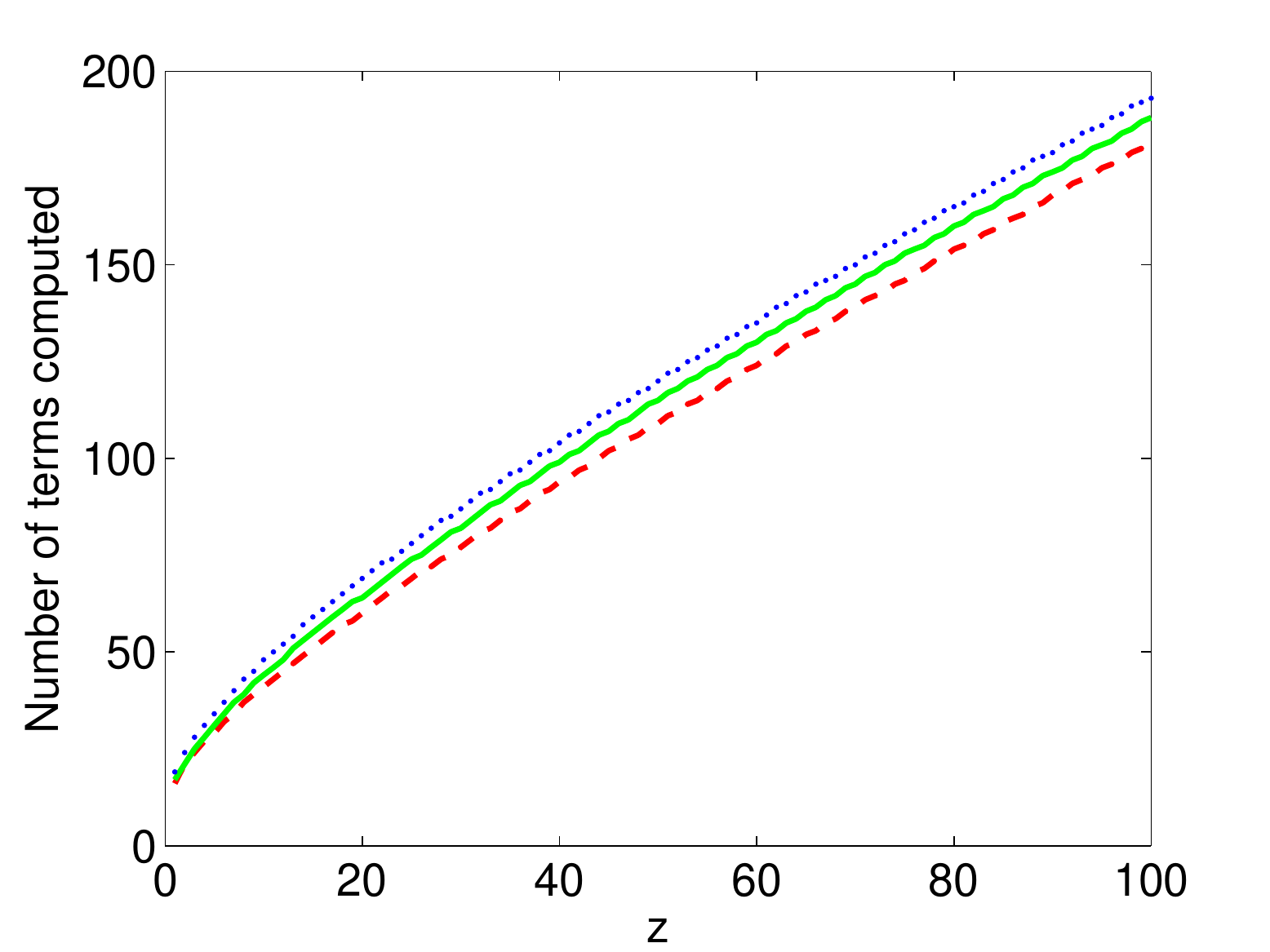}
\caption{(Color online) Number of terms that need to be computed using Taylor series method (a) for evaluating $_{1}F_{1}$ (and hence $\mathbf{M}$) for real $z\in[1,100]$ when $a=2$ and $b=3$ (red, dashed), $a=2+10i$ and $b=10+5i$ (green, solid), and $a=20$ and $b=15$ (blue, dotted) before the stopping criterion (two successive small terms) is satisfied.  The method gives results to 15-digit accuracy.
}
\label{fig:figure3.3}
\end{figure}

By examining the coefficients of the series, we see that $S_N$ converges quickly when $z$ is close to zero or when $b$ has a large real part and a relatively small imaginary part (because $\Gamma(\cdot)^{-1}$ decays rapidly near the positive real axis); again see Table \ref{appbtable1}. However, as $a$ or $z$ becomes large, the coefficients of the Taylor series become large, slowing down the convergence rate (see Fig.~\ref{fig:figure3.3}).  This can be seen in the fourth row in Table \ref{appbtable1}, where the Taylor series is accurate for up to and beyond $\left|z\right|=100$ (although the approximation becomes less accurate as $\left|z\right|\rightarrow\infty$), depending on the parameter values involved.


As demonstrated by cases 13--16 in Table \ref{appbtable1}, a significant issue with computing $\mathbf{M}$ using Taylor series is cancellation.  Although $|a|$ and $|z|$ are the same in each of these cases, the Taylor series method is effective in these cases if $a$ and $z$ have the same signs, whereas it is ineffective if their signs are different. In this situation, we find that the Buchholz polynomial expansion that we describe in Section \ref{sec:1F1_buchholz} is more effective.

One could exploit sequence acceleration methods for computing the series \eqref{1F1defnew}. For example, in the NAG Library routine \url{c06ba} \cite{NAGLibrary}, Shanks' transformation \cite{Shanks} is applied to the sequence using the $\epsilon$-algorithm in Ref.~\cite{Wynn}.


\subsubsection{Writing $S_N$ as a Single Fraction} \label{sec:1F1_singlefraction}

As illustrated by the performance of Taylor series method (a) on cases such as $(a,b,z)=(10^{-8},10^{-12},-10^{-10}+10^{-12}i)$, where machine precision accuracy is not attained despite the apparent relative simplicity of the parameter regime, the methods of Section \ref{sec:1F1_taylor} can often have significant inaccuracies when the parameter values have small modulus, even though one might na\"ively expect the computation to be fairly straightforward. The method that we discuss in this section aims to provide an alternative to those in Section \ref{sec:1F1_taylor}.  It is also based on the basic series definition \eqref{1F1} of the confluent hypergeometric function \cite{NPB1,NPB2}.

The motivation behind this method is that significant round-off error can arise when repeatedly dividing during the computation of individual terms of Taylor series using other methods. Therefore, applying a method that only requires a single division to compute an approximation to $_{1}F_{1}(a;b;z)$ is potentially advantageous.

In this method, we express the series expansion of $_{1}F_{1}(a;b;z)$ as a single fraction rather than as a sum of many fractions. For example, the sum $S_{j}$ of the first $j+1$ terms of the series up to the term in $z^{2}$ can be expressed as
\begin{align*}
	S_{0}=\frac{0+1}{1}\,,~~S_{1}=\frac{b+az}{b}\,,~~S_{2}=\frac{(b+az)(2)(b+1)+a(a+1)z^{2}}{2b(b+1)}\,,
\end{align*}
where the numerators and denominators of $S_{0},S_{1},S_{2},...$ can be calculated using \eqref{nardineq1}. If we take $\alpha_{0}=0$, $\beta_{0}=1$, $\gamma_{0}=1$, and $S_{0}=1$ and define $S_{j}$ to be the $j$th approximation of $_{1}F_{1}(a;b;z)$, we can apply the following algorithm:

\vspace{0.75em}
\underline{\textbf{Method (c):}}
\begin{align}\label{nardineq1} 
	\alpha_{j}&=(\alpha_{j-1}+\beta_{j-1})j(b+j-1)\,, \\
	\nonumber \beta_{j}&=\beta_{j-1}(a+j-1)z\,, \\
	\nonumber \gamma_{j}&=\gamma_{j-1}j(b+j-1)\,, \\
	\nonumber S_{j}&=\frac{\alpha_{j}+\beta_{j}}{\gamma_{j}}\,, \quad{} j=1,2, \ldots ~\,.
\end{align}
Implementing this method generates a sequence of approximations to $_{1}F_{1}(a;b;z)$ (i.e., approximations of $\{S_{j}$ for $j=1,2, \ldots \}$), which we terminate using the same stopping criterion as in Section \ref{sec:1F1_taylor}.

From Table \ref{appbtable1}, one can conclude that the methods of Section \ref{sec:1F1_taylor} generate accurate computations of $_{1}F_{1}(a;b;z)$ more successfully than the methods previously discussed in this section, for a wide range of parameter and variable values. One possible explanation for this is that---especially as the moduli of the parameter values become increasingly large---the terms in the numerator and denominator of $S_{j}$ become very large for a relatively small $j$, so numerical issues such as round-off error and cancellation may become significant when carrying out the computations. Additionally, we often find that either the numerator or the denominator of $S_{j}$ becomes very large after only a few summations.

Nevertheless, the method that we have just described is useful if $\left|b\right|$ is small (especially when $\left|b\right|\lesssim1$), provided $\left|a\right|$ is not too large. However, if one uses the Taylor series methods of Section \ref{sec:1F1_taylor} when $\left|b\right|$ is small, then the round-off error can become costly if many terms are too large. Therefore, a method with a single division can potentially aid accurate computation in this case, as the effect of round-off error is reduced.


\subsection{Buchholz Polynomials} \label{sec:1F1_buchholz}

References \cite{AbadSesma1,AbadSesma2} used Buchholz polynomials to derive expressions for $\mathbf{M}(a;b;z)$, focusing on the regime of real parameters and variable. We find that this yields a much more effective approach than the methods discussed in Sections \ref{sec:1F1_taylor} and \ref{sec:1F1_singlefraction} for the computation of $\mathbf{M}(a;b;z)$ when $\text{sign}(a)=-\text{sign}(z)$.

The function $\mathbf{M}(a;b;z)$ can expanded in terms of \textbf{Buchholz polynomials} $p_{j}(b,z)$ as follows:
\begin{align}\label{Buchholz1} 
\ \mathbf{M}(a;b;z)=e^{z/2}2^{b-1}\sum_{j=0}^{\infty}p_{j}(b,z)\frac{J_{b-1+j}(\sqrt{z\{2b-4a\}})}{(z\{2b-4a\})^{\frac{1}{2}(b-1+j)}}
\end{align}
provided $b\neq2a$, and when $b=2a$  \cite{AbramowitzStegun}
\begin{align*}
 	\mathbf{M}(a;2a;z)=\Gamma\left(a+\frac{1}{2}\right)\Gamma(2a)e^{z}\left(\frac{iz}{4}\right)^{-(a-\frac{1}{2})}J_{a-\frac{1}{2}}\left(\frac{iz}{2}\right),
\end{align*}
where $J_{\nu}(z)$ denotes the Bessel function of the first kind.  The coefficients of (\ref{Buchholz1}) are
\begin{align} \label{Buchholzpoly} 
	p_{j}(b,z)=\frac{(iz)^{j}}{j!}\sum_{s=0}^{\left\lfloor\frac{j}{2}\right\rfloor}\left(\begin{array}{c}
j \\ 2s \\
\end{array}\right)f_{s}(b)g_{j-2s}(z)\,,
\end{align}
where
\begin{align*}
\ f_{0}(b)=1,\quad&f_{s}(b)=-\left(\frac{b}{2}-1\right)\sum_{j=0}^{s-1}\left(\begin{array}{c}
2s-1 \\ 2j \\ \end{array}\right)\frac{4^{s-j}\left|B_{2(s-j)}\right|}{s-j}f_{j}(b)\,,~s=1,2,\ldots~\,, \\
\ g_{0}(z)=1,\quad&g_{s}(z)=-\frac{iz}{4}\sum_{j=0}^{\left\lfloor\frac{s-1}{2}\right\rfloor}\left(\begin{array}{c} s-1 \\ 2j \\ \end{array}\right)\frac{4^{j+1}\left|B_{2(j+1)}\right|}{j+1}g_{s-2j-1}(z)\,,~~s=1,2,\ldots~\,.
\end{align*}
The coefficients $B_{j}$ denote the \textbf{Bernoulli numbers}, which are defined by the generating function
\begin{align*}
	\frac{z}{e^{z}-1}=\sum_{j=0}^{\infty}B_{j}\frac{z^{j}}{j!}\,.
\end{align*}

Another form of equation \eqref{Buchholz1} is as follows for real $a$ (no restrictions are given on the values of $b$ and $z$) \cite{AbadSesma2}:
\begin{align} \label{Buchholz2} 
 	\mathbf{M}(a;b;z)=e^{z/2}2^{b-1}\sum_{j=0}^{\infty}D_{j}z^{j}\frac{J_{b-1+j}(\sqrt{z\{2b-4a\}})}{(z\{2b-4a\})^{\frac{1}{2}(b-1+j)}}\,,
\end{align}
where
\begin{align}
 	\notag D_{0}&=1\,,~~~D_{1}=0\,,~~~D_{2}=\frac{b}{2}\,, \\
 jD_{j}&=(j-2+b)D_{j-2}+(2a-b)D_{j-3}\,,~~j=3,4, \ldots~\,.
\end{align}

The expression in \eqref{Buchholz2} provides an expansion for the confluent hypergeometric function in terms of Bessel functions, which are signficantly easier to compute than $\mathbf{M}$. (See Appendix \ref{sec:appendixc} for more details.)

As illustrated in Table \ref{appbtable1}, this method seems to outperform the Taylor series methods of Section \ref{sec:1F1_taylor} for computation of $\mathbf{M}$ for certain parameter regimes, and this is particular true for many examples in which $\text{sign}(a)\neq\text{sign}(z)$. As illustrated by the first row of Table \ref{appbtable1}, one problematic situation for this method occurs when $b\approx2a$.  This results from the division by powers of $\sqrt{z\{2b-4a\}}$ in \eqref{Buchholz2}. In such cases, one can make use of recurrence relations, which we discuss in Section \ref{sec:1F1_recurrences}, or compute the Taylor series directly.

This method is  valuable for moderate values of $a$ and $z$ (in particular, for $10\lesssim\left|a\right|,\left|z\right|\lesssim100$), especially when the real parts of $a$ and $z$ have opposite signs. The Taylor series methods discussed in Section \ref{sec:1F1_taylor} and the single fraction method of Section \ref{sec:1F1_singlefraction} do not give accurate computations with these cases, but the method described in this section performs very well. (See Table \ref{appbtable1} in Appendix \ref{sec:appendixb}.)  One even obtains a very accurate result when $\left|a\right|$ is large (e.g., when $|a|=500$), with real $z$ of opposite sign. (See case 19 in Appendix \ref{sec:appendixa}, for which we obtained 16-digit accuracy.)  As illustrated in Appendix \ref{sec:appendixb}, however, this method becomes less successful as $\left|z\right|$ becomes larger. Nonetheless, the good performance of this method for a large range of parameter values makes it a convenient element of one's toolbox for computing $\mathbf{M}(a;
b;z)
$.


\subsection{Asymptotic Series for Large $|z|$} \label{sec:1F1_asymptotic}

The methods that we outlined in Sections \ref{sec:1F1_taylor} and \ref{sec:1F1_buchholz} were all ineffective (in general) for large values of $\left|z\right|$.  (We found that these methods ceased to be effective for $\left|z\right|\gtrsim100$, although the threshold was sometimes lower, depending on the precise values of the parameters.) In this section, we aim to address this issue by using large-$\left|z\right|$ asymptotic formulas for computing the confluent hypergeometric function.

One may apply Watson's lemma \cite{Watson} to show that as $\left|z\right|\rightarrow\infty$, the hypergeometric function $\mathbf{M}(a;b;z)$ satisfies \cite{AbramowitzStegun}
\begin{align} \label{asym3} 
\	\mathbf{M}(a;b;z)\sim\frac{e^{z} z^{a-b}}{\Gamma(a)}&\sum_{j=0}^{\infty}\frac{(b-a)_{j}(1-a)_{j}}{j!}z^{-j} \\
\ \notag &+\frac{e^{\pm\pi ia}z^{-a}}{\Gamma(b-a)}\sum_{j=0}^{\infty}\frac{(a)_{j}(1+a-b)_{j}}{j!}(-z)^{-j}\,,
\end{align}
where $a,b\in\mathbb{C}$, $a,b-a\neq\mathbb{Z}^{-}\cup\{0\}$, $-\frac{1}{2}\pi+\delta\leq\pm\arg z\leq\frac{3}{2}\pi-\delta$,\footnote{As in \cite{AbramowitzStegun} where `$\pm$' is stated, the sign `$+$' is used if $-\frac{1}{2}\pi<\arg z<\frac{3}{2}\pi$, with `$-$' used if $-\frac{3}{2}\pi<\arg z<-\frac{1}{2}\pi$. The argument is assigned limits $\pm\left(\frac{3}{2}\pi-\delta\right)$, $\pm\left(-\frac{1}{2}\pi+\delta\right)$ to ensure validity of the expansion over the branch cut present.} and $0<\delta\ll1$.

The expression when $z$ is real is given by \cite{DLMF}
\begin{align*}
 	\mathbf{M}(a;b;z)\sim{}\frac{e^{z}z^{a-b}}{\Gamma(a)}\sum_{j=0}^{\infty}\frac{(b-a)_{j}(1-a)_{j}}{j!}z^{-j}\,.
\end{align*}

We computed these series in {\scshape Matlab} using the same two techniques as for the Taylor series expansion in Section \ref{sec:1F1_taylor}.  In the first method, we compute each term by adding each new term recursively to the previous sum, where we continue including new terms until they become small, using the same stopping criteria as before. In the second method, we find each term iteratively in terms of the previous two terms and then add the new term to the previous sum. As in Section \ref{sec:1F1_taylor}, we label these techniques as methods (a) and (b).  We show our results in Table \ref{appbtable1}.

As with the Taylor series techniques of Section \ref{sec:1F1_taylor}, we obtain similar results using methods (a) and (b). Both methods work well for large $|z|$ and moderate values of the parameters $a$ and $b$.

Because the asymptotic series are expressed in terms of hypergeometric series of the form $_{2}F_{0}$ instead of $_{1}F_{1}$ (or $\mathbf{M}$), there is no longer a $(b)_{j}$ term in the denominator of the terms of the series, so parameter regimes involving $b$ with large modulus can no longer be treated in a straightforward manner.\footnote{We define the hypergeometric function $_{2}F_{0}$ as follows:
\begin{align*}
 	_{2}F_{0}(a,b;-;z)=\sum_{j=0}^{\infty}(a)_{j}(b)_{j}~\frac{z^{j}}{j!}\,.
\end{align*}} The cases that we tested suggest that the methods cope reasonably well when neither $\left|a\right|$ nor $\left|b\right|$ is very large. (As a guide, the computations seem to lose accuracy when $\left|a\right|$ or $\left|b\right|$ is greater than 50, although this depends on the value of $z$. For these cases, one needs to apply recurrence relations; see the discussion in Section \ref{sec:1F1_recurrences} \cite{NPB1,NPB2}.)

References~\cite{OldeDaalhuis1,Olver2,Olver5} discussed the expansion
\begin{align} \label{hyperexp} 
	U(a;b;z)=z^{-a}\sum_{j=0}^{N-1}\frac{(a)_{j}(a-b+1)_{j}}{j!}(-z)^{-j}+R_{N}(a;b;z)\,,
\end{align}
where
\begin{align}
 	R_{N}(a;b;z) &= {}\frac{2\pi(-1)^{N}z^{a-b}}{\Gamma(a)\Gamma(a-b+1)} \\
\ \nonumber \times&{}\left(\sum_{j=0}^{M-1}\frac{(1-a)_{j}(b-a)_{j}}{(-z)^{j}j!}G_{N+2a-b-j}(z)+(1-a)_{M}(b-a)_{M}R_{M,N}(a;b;z)\right)\,, \\
 	\notag G_{\eta}(z) &= {}\frac{e^{z}}{2\pi}\Gamma(\eta)\Gamma(1-\eta,z)\,,
\end{align}
and $\Gamma(\varpi,z)$ is the (upper) incomplete gamma function. Using this expansion, one can compute $\mathbf{M}$ via the expression \eqref{bfMU}, though we find the problem to be more tractable when using other approaches described in this section.

As stated in Ref.~\cite{Olver5}, one can take $\delta$ to be an arbitrary small parameter and consider $\left|z\right|\rightarrow\infty$ with $a$, $b$, and $m$ fixed to obtain
\begin{align} \label{rmn} 
	R_{M,N}(a;b;z)=\left\{\begin{array}{ll}
O(e^{-\left|z\right|}z^{-M})\,, & \text{if}~\left|\arg z\right|<\pi\,,\\
O(e^{z}z^{-M})\,, & \text{if}~\pi\leq\left|\arg z\right|\leq\frac{5}{2}\pi-\delta\,.\\
\end{array}~\right.
\end{align}
This gives a \emph{uniform, exponentially improved asymptotic expansion} \cite{Olver2}.  A second expansion with exponentially-improved accuracy is derived in Ref.~\cite{OldeDaalhuis1}.  This was in turn extended in Ref.~\cite{OldeDaalhuisOlver} to a \textbf{hyperasymptotic expansion}, which offers another possibility for computing $\mathbf{M}$ for large $|z|$.

There has also been research on uniform asymptotic expansions of $\mathbf{M}(a;b;z)$ for large values of the parameters $a$ and $b$.  For example, reference \cite{Temme6} derived an expansion for large (real and positive) values of $b$ and $z$ in terms of parabolic cylinder functions. (We discuss methods for computing parabolic cylinder functions in Appendix \ref{sec:appendixc}.)  Reference \cite{Temme7} gave an expansion for $U(a;b;z)$ for large values of $a$ (which must be real and positive) and $b\leq1$ in terms of modified Bessel functions (which are themselves expressible in terms of Bessel functions \cite{DLMF}). To our knowledge, no similar uniform asymptotic expansions have been derived for complex parameter or variable values.


\subsection{Quadrature Methods} \label{sec:1F1_quadrature}

Thus far, we have only considered methods based on series. In this section, we discuss a class of methods for computing $\mathbf{M}(a;b;z)$ using its integral representation for ${\rm Re}(b)>{\rm Re}(a)>0$.

When ${\rm Re}(b)>{\rm Re}(a)>0$, the function $\mathbf{M}(a;b;z)$ has the integral representation \cite{AbramowitzStegun}
\begin{align} \label{1F1int} 
\	\mathbf{M}(a;b;z)=\frac{1}{\Gamma(a)\Gamma(b-a)}\int_{0}^{1}e^{zt}(1-t)^{b-a-1}t^{a-1}{\rm d}t\,.
\end{align}
Applying the transformation $t\mapsto\frac{1}{2}\widetilde{t}+\frac{1}{2}$ and using Jacobi parameters $\widetilde{\alpha}=b-a-1$ and $\widetilde{\beta}=a-1$ yields \cite{Gautschi3}
\begin{align*} \label{1F1int10*} 
	\mathbf{M}(a;b;z) &= \frac{1}{2^{b-1}\Gamma(a)\Gamma(b-a)}\int_{-1}^{1}e^{z\left(\frac{1}{2}\widetilde{t}+\frac{1}{2}\right)}(1-\widetilde{t})^{b-a-1}(1+\widetilde{t})^{a-1}{\rm d}\widetilde{t} \\
\ &=\frac{1}{\Gamma(a)\Gamma(b-a)}\left[\frac{e^{z/2}}{2^{b-1}}\sum_{j=1}^{N_{\mbox{\scriptsize mesh}}}w_{j}^{GJ}e^{zt_{j}^{GJ}/2}+E_{N_{\mbox{\scriptsize mesh}}}(a;b;z)\right]\,,
\end{align*}
where $t_{j}^{GJ}$ and $w_{j}^{GJ}$ are the Gauss-Jacobi nodes and weights on $[-1,1]$, and $N_{\mbox{\scriptsize mesh}}$ is the number of mesh points.

This method is known as \textbf{Gauss-Jacobi quadrature}. The error $E_{N_{\mbox{\scriptsize mesh}}}$ for this method can be controlled by $N_{\mbox{\scriptsize mesh}}$, which one can see for real $a$ and $b$ using the following relation \cite{Gautschi3}:
\begin{align} 
	N_{\mbox{\scriptsize mesh}}{}\geq\frac{e\left|z\right|}{8}~\widetilde{t}\left(\frac{4}{e\left|z\right|}\left[x_{+}+(3-2b)\log2+\log\left(\frac{1}{E_{N_{\mbox{\scriptsize mesh}}}}\right)\right]\right)\,,
\end{align}
where $x={\rm Re}(z)$, and $x_{+}$ is equal to $x$ if $x\geq0$ and $0$ if it is negative. The quantity $\widetilde{t}$ denotes the inverse of the function $s=\widetilde{t}\log \widetilde{t}$. Low-order approximations are given in Ref.~\cite{Gautschi1} for different real values of $s$.

Gauss--Jacobi quadrature is a natural choice because of the form of the integrand in \eqref{1F1int} and the fact that the integrand blows up at the end points of the integral. To implement Gauss--Jacobi quadrature, we use the Golub-Welsch algorithm \cite{GolubWelsch}, though of course other methods are possible (e.g., the Glaser-Liu-Rokhlin algorithm \cite{GlaserLiuRokhlin}, or the Hale-Townsend algorithm \cite{HaleTownsendQuad}). As illustrated in Table \ref{appbtable1}, we find that using Gauss-Jacobi quadrature is effective for a large range of values of $\left|z\right|$, provided $z$ does not have an imaginary part with magnitude greater than roughly 100.  Additionally, a problem arises when either $\left|a\right|$ or $\left|b\right|$ becomes fairly large. A problem also occurs if $|z|$ is large, as the integral under consideration is then highly oscillatory and/or stiff.  

For small values of $\left|a\right|$ and $\left|b\right|$ (up to about 30--40), the method of Gauss-Jacobi quadrature is extremely useful for evaluating the confluent hypergeometric function when ${\rm Re}(b)>{\rm Re}(a)>0$, as the methods that we described in previous sections are not always reliable in this parameter regime. Due to the computational expense of generating the quadrature nodes and weights, which is far greater than the cost of implementing many of the series expansion methods previously discussed, we would only recommend a quadrature approach in parameter regimes where series methods are likely to struggle.

Gauss-Jacobi quadrature is the most effective quadrature method among those that we investigated.  It is, however, useful to mention others briefly.  For example, $U(a;b;z)$ can be computed using a trapezoidal-rule method \cite{AllasiaBesenghi}.  Additionally, Ref.~\cite{Nieuwveldt} described how to use a contour integral method to evaluate $\mathbf{M}(a;b;z)$, Refs.~\cite{SchmelzerTrefethen,TWS,Weideman,WeidemanTrefethen} implemented a Talbot contour integral method, and Ref.~\cite{Michel} used direct integration via continued fractions.


\subsection{Recurrence Relations} \label{sec:1F1_recurrences}

When computing the confluent hypergeometric function, a method's robustness---i.e., the range of parameter and variable values for which it is effective---is often greatly reduced by its poor performance as $\left|{\rm Re}(a)\right|$ or $\left|{\rm Re}(b)\right|$ increases.  In this section, we detail recurrence relation techniques, which can reduce the problem of computation with large parameter values to a simpler problem of computing $\mathbf{M}(a;b;z)$ with values of $|{\rm Re}(a)|$ and $|{\rm Re}(b)|$ much closer to $0$. One of the methods discussed in Sections \ref{sec:1F1_taylor}--\ref{sec:1F1_quadrature} can subsequently be applied to solve the simpler problem, usually with much greater success than with a direct computation.

The function $M(a;b;z)$ satisfies the following recurrence relations \cite{GST2,GST3,SeguraTemme}:
\begin{align}
	 \nonumber (b-a)~M(a-1;b;z)+(2a-b+z)~M(a;b;z)-a~M(a+1;b;z)&=0\,, \\
	 \nonumber b(b-1) ~M(a;b-1;z)+b(1-b-z)~M(a;b;z)+z(b-a) ~M(a;b+1;z)&=0\,, \\ 
	 \label{1F1rec3} b(b-1)~M(a-1;b-1;z)-b(b-z-1)~M(a;b;z)-a z~M(a+1;b+1;z)&=0\,,
\end{align}
from which $\mathbf{M}$ can be computed. We denote these recurrence relationships using the standard notation $(+0)$, $(0+)$, and $(++)$ \cite{GST2,SeguraTemme}. The $+$  indicates which parameters ($a$, $b$, or both) are increasing.  A na\"ive way of using these relationships is via direct application.  However, for certain parameter regimes, it is instead appropriate to utilize the theory of \textbf{minimal solutions} of a recurrence relationship.

A solution $f_{n}$ of a recurrence relation
\begin{align} \label{genrec} 
	y_{n+1}+b_{n}y_{n}+a_{n}y_{n-1}=0
\end{align}
is said to be a \textbf{minimal solution} if there exists a linearly independent solution $g_{n}$ (called a \textbf{dominant solution}) such that
\begin{align*}
	\lim_{n\rightarrow\infty}\frac{f_{n}}{g_{n}}=0\,.
\end{align*}
We use the following theorem, discussed in \cite{GST2,GST3,SeguraTemme}, in our subsequent investigation of recurrence relations.

\vspace{0.75em}
\begin{Thm} \label{Poincare}
(\textbf{Poincar\'{e}'s Theorem}) ~Consider the recurrence relation \eqref{genrec}, where $\lim_{n\rightarrow+\infty}b_{n}=b_{\infty}$ and $\lim_{n\rightarrow+\infty}a_{n}=a_{\infty}$. Denote the zeros of the equation $t^{2}+b_{\infty}t+a_{\infty}=0$ by $t_{1}$ and $t_{2}$. If $\left|t_{1}\right|\neq\left|t_{2}\right|$, then the recurrence relation \eqref{genrec} has two linearly independent solutions $f_{n}$, $g_{n}$ such that
\begin{align*}
	\lim_{n\rightarrow+\infty}\frac{f_{n}}{f_{n-1}}=t_{1}\,,\quad\lim_{n\rightarrow+\infty}\frac{g_{n}}{g_{n-1}}=t_{2}\,.
\end{align*}
If $\left|t_{1}\right|=\left|t_{2}\right|$, then
\begin{align*}
 	\limsup_{n\rightarrow+\infty}\left|y_{n}\right|^{1/n}=\left|t_{1}\right|
\end{align*}
for any nontrivial solution $y_{n}$ of \eqref{genrec}.

Finally, when $\left|t_{1}\right|\neq\left|t_{2}\right|$, the solution whose ratio of consecutive terms tends to the root of smaller modulus is always the minimal solution. $\Box$
\end{Thm}
\vspace{0.75em}

We now consider the recurrence relations \eqref{1F1rec3}, which may be modified to relations of the form \eqref{genrec}, with
\begin{align}
\ \label{anbn1F1} \underline{(+0):}\quad\quad&a_{n}=\frac{a-b+n}{a+n},\quad\quad~~b_{n}=-\frac{2a-b+z+2n}{a+n}, \\
\ \nonumber \underline{(0+):}\quad\quad&a_{n}=\frac{b-a+n-1}{z},\quad{}b_{n}=\frac{1-b-n-z}{z}, \\
\ \nonumber \underline{(++):}\quad\quad&a_{n}=-\frac{1}{(a+n)z},\quad\quad~b_{n}=\frac{b-z-1+n}{(a+n)z}.
\end{align}
The solutions to the recursion relations \eqref{anbn1F1} and their ratios as $n\rightarrow+\infty$  (for real $z$) become \cite{GST3}
\begin{align} \label{abovefact1} 
	&f_{n}=\Gamma(1+a+n-b)U(a+n;b;z)\,,~g_{n}=M(a+n;b;z)\,,~\frac{f_{n}}{g_{n}}\sim e^{-4\sqrt{nz}}\,, \\
	 \nonumber &f_{n}=\Gamma(b+n-a)\mathbf{M}(a;b+n;z)\,,~g_{n}=U(a;b+n;z)\,,~\frac{f_{n}}{f_{n-1}}\sim1,~\frac{g_{n+1}}{g_{n}}\sim\frac{n}{z}\,, \\
	 \nonumber &f_{n}=\mathbf{M}(a+n;b+n;z)\,,~g_{n}=(-1)^{n}U(a+n;b+n;z)\,,~\frac{f_{n}}{f_{n-1}}\sim\frac{1}{n}\,,~\frac{g_{n}}{g_{n-1}}\sim-\frac{1}{z}\,.
\end{align}

Using \eqref{abovefact1}, the definition of a minimal solution, and Poincar\'{e}'s Theorem, one can deduce that the minimal solutions of the relations $(+0)$, $(0+)$, and $(++)$ are, respectively,
\begin{align} \label{1F1minsols} 
	\Gamma(1+a+n-b)U(a+n;b;z)\,,~\Gamma(b+n-a)\mathbf{M}(a;b+n;z)\,,~\mathbf{M}(a+n;b+n;z)\,.~~~~~~~~~~
\end{align}
When applying the $(+0)$ recursion, however, the solution is only minimal in $\mathbb{C}\backslash\mathbb{R}$ \cite{SeguraTemme}, making this a less convenient approach in the large $|a|$ case.

\subsubsection{Miller's Algorithm}

Suppose that we seek a solution of the three-term recurrence relation \eqref{genrec}. One method to compute numerical approximations of $\widetilde{f}_{n}$ $(n=0,\ldots,k)$ for a minimal solution $f_{n}$ (if the recurrence admits such a solution\footnote{Where a minimal solution exists, it corresponds to $f_{n}$ in the notation of this section.}) is to apply \textbf{Miller's algorithm} \cite{GST3}. 
The solution $f_k \over f_N$ satisfies the following system:
 	$$\begin{pmatrix}
		a_1 & b_1 & 1 \cr
		&a_2 & b_2 &  1 \cr
		&& \ddots & \ddots & \ddots \cr
		&&& a_{N-1} & b_{N-1} & 1 \cr
		&&&&1\cr
		&&&&&1
	\end{pmatrix}  \begin{pmatrix}
		{f_0 / f_N} \cr
		{f_1 / f_N } \cr 
		\vdots \cr
		{f_{N-1} / f_N} \cr
		{1} \cr
		{f_{N+1} / f_N }
	\end{pmatrix} = \begin{pmatrix} 0 \cr 0 \cr \vdots \cr 0 \cr 1 \cr f_{N+1}/f_N \end{pmatrix}.$$
Miller's algorithm involves observing that $f_{N+1}/f_N \rightarrow 0$ as $N \rightarrow \infty$.  Thus we can approximate $f_k/f_N$ by $y_k^N$ which solves
 	$$\begin{pmatrix}
		a_1 & b_1 & 1 \cr
		&a_2 & b_2 &  1 \cr
		&& \ddots & \ddots & \ddots \cr
		&&& a_{N-1} & b_{N-1} & 1 \cr
		&&&&1\cr
		&&&&&1
	\end{pmatrix}  \begin{pmatrix}
		y_0^N \cr
		y_1^N \cr 
		\vdots \cr
		y_{N-1}^N \cr
		y_N^N \cr
		y_{N+1}^N
	\end{pmatrix} = \begin{pmatrix} 0 \cr 0 \cr \vdots \cr 0 \cr 1 \cr 0 \end{pmatrix}.$$
Then we obtain the approximation 
	$f_k \approx {f_0 \over y_N^N} y_k^N.$
The issue is choosing $N$ so that this approximation is within a given tolerance.    The standard procedure is to continually increase $N$ until  the desired $f_k$ has changed less than a proscribed tolerance, i.e., to verify:
	$$\left|{f_0 \over y_{N+1}^{N+1}} y_k^{N+1}  - {f_0 \over y_{N}^{N}} y_k^{N}\right| < {\tt tol}.$$

\subsubsection{Olver's Algorithm}

 \def\vc#1{{\mathbf #1}}

 Another method for the computation of the minimal solution of \eqref{genrec} is to apply Olver's algorithm \cite{Olver3,Wimp}, which avoids the issue of choosing $N$.   Olver's algorithm treats the recurrence relationship as a two-point boundary value problem: for any $N$,  the minimal solution satisfies
 	$$\begin{pmatrix}
		1 \cr
		a_1 & b_1 &  1 \cr
		& a_2 & b_2 & 1 \cr
		&& \ddots & \ddots & \ddots \cr
		&&& a_{N-1} & b_{N-1} & 1 \cr
		&&&&&1
	\end{pmatrix}  \begin{pmatrix}
		f_0 \cr
		{f_1 } \cr 
		f_2 \cr
		\vdots \cr
		{f_{N-1}} \cr
		{f_N }
	\end{pmatrix} = \begin{pmatrix} 1 \cr 0 \cr 0 \cr \vdots \cr 0 \cr {f_N } \end{pmatrix},$$
 %
Now consider solving this linear system (ignoring for now that we do not know ${f_N }$) using Gaussian elimination.  The key observation is that forward elimination (without pivoting) is independent of $N$, thus we can perform it {\it adaptively}, and only performing back substitution when a a convergence criteria is satisfied.  This convergence criteria can be tested in $O(1)$ operations, {\it without} performing the back substitution.


Explicitly, define recursively
	\begin{align*}
		p_1 = 1,\qquad p_2 &= -b_1, \qquad p_{k+1} = -a_{k} p_{k-1} - b_{k} p_k, \cr
		r_1 &= a_1,\qquad r_{k+1} = a_{k+1} r_k.
	\end{align*}
%

After $N+1$ iterations of Gaussian elimination we have the system
 	$$\begin{pmatrix}
		1 \cr
		 & p_2 &  -p_1 \cr
		&  & p_3 & -p_2 \cr
		&&  & \ddots & \ddots \cr
		&&&  & p_{N-1} & -p_{N-2} \cr				
		&&&&  & p_{N} & -p_{N-1} \cr	
		&&&&&  & p_{N+1} & -p_{N} \cr				
		&&&&&& a_{N+2} & b_{N+2} & 1 \cr
		&&&&&&& \ddots & \ddots & \ddots
	\end{pmatrix}   \begin{pmatrix}
		f_0 \cr
		{f_1 } \cr 
		f_2 \cr
		\vdots \cr
		f_{N-1} \cr
		{f_{N}} \cr
		{f_{N+1} } \cr
		f_{N+2} \cr
		\vdots
	\end{pmatrix}  = \begin{pmatrix} 1 \cr r_1 \cr r_2 \cr \vdots \cr r_{N-1} \cr r_{N} \cr r_{N+1} \cr 0 \cr \vdots \end{pmatrix}.$$
 (Note: this is a slight variant of forward elimination that introduces a zero in the $k$th column of the matrix by multiplying  the $(k+1)$th row of the system by $-p_k$, then adding $a_k$ times the $k$th row to the $(k+1)$th row.) 
 
For any fixed $N$, we can truncate the system so it becomes upper triangular
$$\begin{pmatrix}
		1 \cr
		 & p_2 &  -p_1 \cr
		&  & p_3 & -p_2 \cr
		&&  & \ddots & \ddots \cr
		&&&  & p_{N-1} & -p_{N-2} \cr				
		&&&&  & p_{N} 
	\end{pmatrix}   \begin{pmatrix}
		f_1^N \cr
		{f_2^N } \cr 
		f_3^N \cr
		\vdots \cr
		f_{N-1}^N \cr
		{f_{N}^N} \cr
	\end{pmatrix}	= \begin{pmatrix} 1 \cr r_2 \cr r_3 \cr \vdots \cr r_{N-1} \cr r_{N} \end{pmatrix}$$
and then $f_k^N$ are obtainable by back substitution.  

 Now a sensible way for choosing $N$ is such that $|{f_k^{N+1} - f_k^N}|$ is below a given tolerance.  The important observation in \cite{Olver3} is that we do not have to perform back substitution to determine this.  Indeed, we have
\begin{align*}
	&\begin{pmatrix}
		1 \cr
		 & p_2 &  -p_1 \cr
		&  & p_3 & -p_2 \cr
		&&  & \ddots & \ddots \cr
		&&&  & p_{N-1} & -p_{N-2} \cr				
		&&&&  & p_{N} & -p_{N-1} \cr
		&&&&&& p_{N+1} 
	\end{pmatrix}  \left( \begin{pmatrix}
		f_1^{N+1} \cr
		{f_2^{N+1} } \cr 
		f_3^{N+1} \cr
		\vdots \cr
		f_{N-1}^{N+1} \cr
		{f_{N}^{N+1}} \cr
		{f_{N + 1}^{N+1}} 		
	\end{pmatrix} -  \begin{pmatrix}
		f_1^{N} \cr
		{f_2^{N} } \cr 
		f_3^{N} \cr
		\vdots \cr
		f_{N-1}^{N} \cr
		{f_{N}^{N}} \cr
		0
	\end{pmatrix}\right)	\cr
	\qquad & = \begin{pmatrix} 1 \cr r_2 \cr r_3 \cr \vdots \cr r_{N-1} \cr r_{N} \cr r_{N+1} \end{pmatrix} -  \begin{pmatrix} 1 \cr r_2 \cr r_3 \cr \vdots \cr r_{N-1} \cr r_{N} \cr 0\end{pmatrix} = r_{N+1} \vc e_{N+1}.
\end{align*}
Thus the backward error is obtained by applying back-substitution to $ \vc e_{N+1}$, which results in precisely the vector
	$$\left[0,{p_1 \over p_N p_{N+1}}, {p_2 \over p_N p_{N+1}},\ldots,{p_{N-2} \over p_N p_{N+1}}, {p_{N-1} \over p_N p_{N+1}} ,{1 \over p_{N+1}}\right]^\top.$$
 	Thus we know convergence has occurred when
	$$\left|{r_{N+1}  \over p_N p_{N+1}}\right| \max_{k=1,\ldots,N-1} |p_k|  < {\tt tol}.$$
This takes only $O(1)$ operations to determine, versus the $O(N)$ operations of back-substitution.  
 
%

\subsubsection{Application to Hypergeometric Functions}

The technique of using recurrence relations can be applied to compute $\mathbf M$ effectively for parameter regimes for which previously described methods do not perform well. Additionally, the ideas discussed in this section can be extended to compute recurrence relations with large $\left|{\rm Re}(a)\right|$ \cite{GST3}. As discussed in Section \ref{sec:1F1_summary}, one can also apply these relations to successfully rephrase the problem of computing the confluent hypergeometric function as one involving larger values of parameters. 

Note, however, that the minimal solutions \eqref{1F1minsols} involve Gamma functions, so the effectiveness of the aforementioned technique is potentially restricted by issues with computing the Gamma function of a variable with large modulus, unless great care is taken to counter this possibility when selecting the recurrence relation used. For example, in double-precision computations in the software package {\scshape Matlab}, the number $\Gamma(172)$ is stated as being infinite, even though it is actually finite.  Ideally, one would apply the technique of using recurrence relations to software in which Gamma functions of modulus $\lesssim171$ can be computed effectively.  An alternative, which can be applied to any other method which requires the computation of a Gamma function of large modulus, is to instead calculate the logarithm of each of the Gamma functions involved (as well as that of the hypergeometric function that we seek to compute) and apply the relation $\log(z_{1})+\log(z_{2})=\log(z_{1}z_{2})$.


\subsection{Other Methods for Computing $\mathbf{M}$} \label{sec:1F1_othermethods}

A class of methods for computing $\mathbf{M}(a;b;z)$ or $_{1}F_{1}(a;b;z)$ that can be useful for small $\left|z\right|$ is solving the \textbf{confluent hypergeometric differential equation} \eqref{1f1de} numerically. To do this, we express \eqref{1f1de} as a system of differential equations by writing $f'=g$, $g'=-\frac{1}{z}[(b-z)g-af]$, with $f(0)=1$, $g(0)=\frac{a}{b}$, using the Taylor series expansion \eqref{1F1} to establish $g(0)=f'(0)$. We would then solve this system of first order ordinary differential equations numerically for $f(z)$ and $g(z)$, with our approximation of $_{1}F_{1}(a;b;z)$ being the numerical solution for $f(z)$. The code \url{d02bj} from the NAG Library, which contains an error-controlled adaptive Runge-Kutta method \cite{NAGLibrary}, is an effective means to do this.  Because of issues with computation time, however, this is only a valuable method for small $|z|$, and we would certainly not recommend this approach when a (much cheaper) series method generates accurate solutions.

Many other methods have been developed to compute $\mathbf{M}(a;b;z)$ and $_{1}F_{1}(a;b;z)$. We mention several of them in passing; see Ref.~\cite{Pearson} for additional details.  Reference \cite{Muller} gives two additional series expansions for $_{1}F_{1}(a;b;z)$; one expansion is in terms of beta random variables, and the other is in terms of the (lower) incomplete gamma function.  Reference \cite{LopezPagola} gives an asymptotic expansion for large $\left|b\right|$ and $\left|z\right|$, and Ref.~\cite{OldeDaalhuis1} includes the derivation of a second hyperasymptotic expansion. The integral representation \eqref{1F1int} can also be computed using other methods, such as splitting the integral, Romberg integration, and adaptive quadrature (although the latter two require there to be no blow up of the integrand on either side of the region of integration).  Additional computational techniques include Pad\'{e} or rational approximations \cite{Luke1,Luke2,Luke5,PTVF}, a continued fraction expression \cite{
Wolfram1} (techniques for solving such expressions are described in Ref.~\cite{AbramowitzStegun}), a Chebyshev expansion for $_{1}F_{1}$ \cite{Luke5}, or other expansions in terms of Bessel functions \cite{Luke3}.

We note that considerable research has investigated how to compute the \textbf{Whittaker function}, which is closely related to the confluent hypergeometric function.  Techniques from such inquiries could perhaps be adapted for the computation of $\mathbf{M}$. Details of methods for computing the Whittaker function can be found in Refs.~\cite{BMOF, Dunster2, Lopez, NobleThompson, Olver1, Olver4, Olver6, DLMF}.


\subsection{Summary and Discussion} \label{sec:1F1_summary}

The methods that we have examined and implemented\footnote{We have made {\scshape Matlab} code available at \cite{PearsonCodeLink}.} for the computation of $\mathbf{M}(a;b;z)$ include the series methods in Sections \ref{sec:1F1_taylor}, \ref{sec:1F1_buchholz}, and \ref{sec:1F1_asymptotic}; the use of quadrature in Section \ref{sec:1F1_quadrature}; recurrence relations in Section \ref{sec:1F1_recurrences}); and various additional methods that are described in Section \ref{sec:1F1_othermethods}.

The series methods seemed to generate the most accurate results, and they had very fast computation times in comparison to the inbuilt {\scshape Matlab} function \url{hypergeom}. Additionally, for values of $\left|a\right|$ and $\left|z\right|$ less than about 50 and $\left|b\right|\gtrsim1$, the Taylor series methods described in Section \ref{sec:1F1_taylor} and the method of expressing $_{1}F_{1}$ as a single fraction in Section \ref{sec:1F1_singlefraction} typically seemed to be sufficiently robust (with the restriction for the single-fraction method that $\left|b\right|\lesssim1$). The method of Gauss-Jacobi quadrature is also effective for `moderate' values of $\left|a\right|$, $\left|b\right|$, and $\left|z\right|$.

One parameter regime in which the Taylor series and single-fraction methods both fail when one might expect them to succeed is when $\left|a\right|$ and $\left|z\right|$ are roughly between 10 and 60 with real parts of opposite signs (the problem arises from cancellation). In this case, the method involving Buchholz polynomials discussed in Section \ref{sec:1F1_buchholz} can be very useful.  One can do effective computations for large $\left|z\right|$---another important case---by applying the asymptotic expansions of Section \ref{sec:1F1_asymptotic}.  In Section \ref{sec:1F1_asymptotic}, we also discussed hyperasymptotics, which have exponentially-improved accuracy.  Their use can provide a viable alternative for the computation of $\mathbf{M}$ for large $\left|z\right|$ in cases in which software with high precision and programs to compute the incomplete gamma function in the entire complex plane are available \cite{Olver2,Olver5}.

If $\left|{\rm Re}(a)\right|\gtrsim50$ or $\left|{\rm Re}(b)\right|\gtrsim50$, one can use the recurrence relations that we detailed in Section \ref{sec:1F1_recurrences} in order to compute hypergeometric functions with parameter values that have real parts of smaller absolute value. One then completes the computation in this new parameter regime using one of the methods discussed above.

When the parameters and variable are all real, the main problem that one faces is cancellation. For example, when computing a series representation of $_{1}F_{1}$ (or $\mathbf{M}$) in which positive and negative terms alternate, we find that the resulting error properties are likely to be poor. For instance, consider the computation of $_{1}F_{1}(50;20;-100)$ or $\mathbf{M}(50;20;-100)$. The sum of all of the positive Taylor series terms for $_{1}F_{1}(50;20;-100)$ is the same as the absolute value of the sum of all of the negative terms (roughly $1.7\times10^{62}$) to 16 digits, whereas the value of $_{1}F_{1}(50;20;-100)$ is roughly $1.4\times10^{-40}$. This problem may be circumvented by applying forward recurrences so that the values of $b$ in the eventual direct computation of $\mathbf{M}$ are (positive and) large enough that the terms cancelling each other no longer have as large a magnitude. We have thereby devised an effective strategy for computing $\mathbf{M}$ for real parameters and variables. We 
use the result \cite{DLMF}
\begin{align} \label{importanttrans1F1} 
	_{1}F_{1}(a;b;z)=e^{z}~_{1}F_{1}(b-a;b;-z)\quad\Leftrightarrow\quad\mathbf{M}(a;b;z)=e^{z}~\mathbf{M}(b-a;b;-z)
\end{align}
and show our aggregate strategy in Table \ref{1f1recommendationsreal}.

\begin{table}[h]
\renewcommand{\arraystretch}{1.15}
\begin{center}
\begin{footnotesize}
\begin{tabular}{|@{}c@{}|@{}c@{}|@{}l@{}|@{}c@{}|}
\hline
~Case~ & ~Regions for $a$, $b$, $z$~ & ~Recommended method(s) & ~Relevant sections and references~ \\ \hline \hline
I(A) & $a,b\geq0$, $z\geq0$ & ~Taylor series methods; & \ref{sec:1F1_taylor} \cite{AbramowitzStegun,Muller} \\
 & & ~Single-fraction method; & \ref{sec:1F1_singlefraction} \cite{NPB1,NPB2} \\
 & & ~Buchholz polynomial method; & \ref{sec:1F1_buchholz} \cite{AbadSesma1,AbadSesma2} \\
 & & ~(Hyper)asymptotics; & \ref{sec:1F1_asymptotic} \cite{AbramowitzStegun,Berry,OldeDaalhuis1,OldeDaalhuisOlver,Olver2,Olver5} \\
 & & ~Gauss-Jacobi quadrature, & \ref{sec:1F1_quadrature} \cite{Gautschi3} \\
 & & ~if ${\rm Re}(b)>{\rm Re}(a)>0$; & \\ \hline
I(B) & $a,b\geq0$, $z<0$ & ~Forward recursion to combat cancellation, & \ref{sec:1F1_recurrences} \cite{GST3,SeguraTemme} \\
 & & ~then Taylor series & \ref{sec:1F1_taylor} \cite{AbramowitzStegun,Muller} \\
 & & ~or single fraction method & \ref{sec:1F1_singlefraction} \cite{NPB1,NPB2} \\ \hline
II & $a<0$, $b\geq0$ & ~Transformation \eqref{importanttrans1F1}, & \ref{sec:1F1_summary} \cite{DLMF} \\
 & & ~then same as for Case I & \\ \hline
III & $a,b<0$ & ~Forward recursion to combat cancellation, & \ref{sec:1F1_recurrences} \cite{GST3,SeguraTemme} \\
 & & ~then Taylor series & \ref{sec:1F1_taylor} \cite{AbramowitzStegun,Muller} \\
 & & ~or single-fraction method & \ref{sec:1F1_singlefraction} \cite{NPB1,NPB2} \\ \hline
IV & $a\geq0$, $b<0$ & ~Recurrence relations; & \ref{sec:1F1_recurrences} \cite{GST3,SeguraTemme} \\
 & & ~Transformation \eqref{importanttrans1F1}, & \ref{sec:1F1_summary} \cite{DLMF} \\
 & & ~then same as for Case III & \\ \hline
\end{tabular}
\caption{Recommendations for which methods to use for computing the confluent hypergeometric function when the parameters and the variable are real.}
\label{1f1recommendationsreal}
\end{footnotesize}
\end{center}
\end{table}

\begin{figure}
\centering
\includegraphics[bb=75 440 850 780,width=22cm]{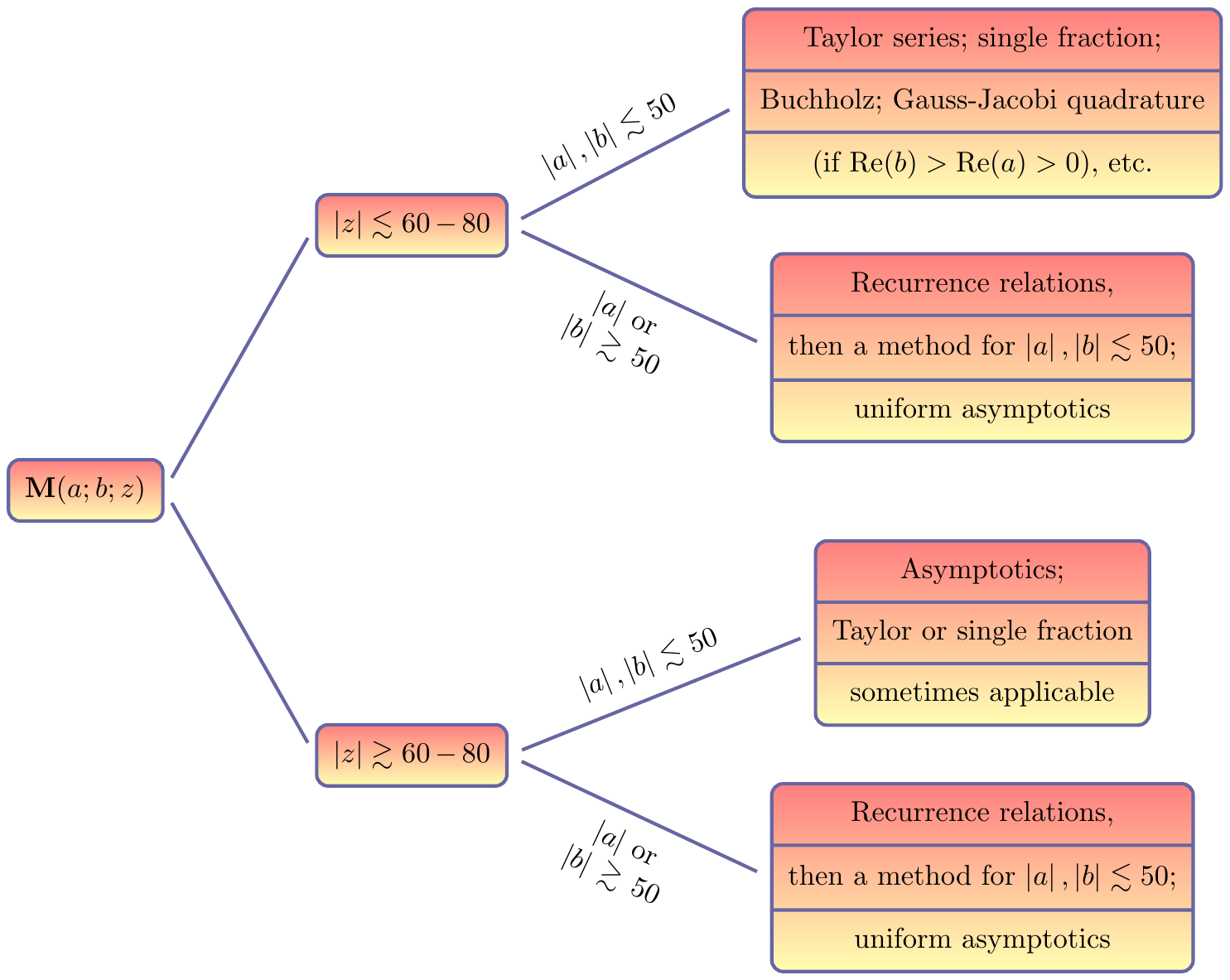}
\caption{Summary of recommended methods for computing $_{1}F_{1}$ for different regions of $a$, $b$, and $z$.
}
\label{fig:summary1F1pic}
\end{figure}

Extending the parameters and variable into the complex plane makes the problem 
of computing the confluent hypergeometric function
much more complicated.  We summarize our recommendations in Fig.~\ref{fig:summary1F1pic}, which shows which of the methods that we have discussed should be used in which regimes.  Details of other software for the computation of the confluent hypergeometric function are discussed in Refs.~\cite{Forrey, Hsu, LozierOlver, MichelStoitsov, Moshier, Temme4}.

\section{Computation of the Gauss Hypergeometric Function $\mathbf{F}$} \label{sec:2F1}

In this section, we discuss methods to compute the Gauss hypergeometric function $\mathbf{F}(a,b;c;z)$ accurately and quickly, and we provide recommendations for the most effective methods for each parameter and variable regime. We implement ideas of a similar nature to those that we discussed for the confluent hypergeometric function---such as those in Sections \ref{sec:2F1_taylor} and \ref{sec:2F1_quadrature}---and we also consider additional methods that are only applicable to $\mathbf{F}$, such as specific transformations and analytic continuation formulas. We then discuss additional methods in Section \ref{sec:2F1_othermethods} and summarize the desired computational strategy in Section \ref{sec:2F1_summary}.


\subsection{Properties of $\mathbf{F}$} \label{sec:2F1_properties}

As we indicated in Section \ref{sec:background}, the \textbf{Gauss hypergeometric function} $\mathbf{F}(a,b;c;z)$ is defined as the series \eqref{2F1} when $z$ is in the radius of convergence (i.e., when $\left|z\right|<1$). This series is defined for any $a\in\mathbb{C}$, $b\in\mathbb{C}$, and $c\in\mathbb{C}$. For $z$ outside of this range, $\mathbf{F}(a,b;c;z)$ is defined by analytic continuation, which we detail in Section \ref{sec:2F1_analyticcont}. This makes it possible to consider the computation of $\mathbf{F}$ for any $z\in\mathbb{C}$.

Note that $\mathbf{F}(a,b;c;0)=\frac{1}{\Gamma(c)}$ for any $a,b \in \mathbb{C}$ and $c\notin\mathbb{Z}^{-}\cup\{0\}$. If $c=n \in\mathbb{Z}^{-}\cup\{0\}$, then \eqref{2F1} is given by a polynomial in $z$ of degree $-n$.  On the unit disk $\left|z\right|=1$, the series in \eqref{2F1} converges absolutely when ${\rm Re}(c-a-b)>0$ \cite{AbramowitzStegun}, and it attains the value $\frac{\Gamma(c)\Gamma(c-a-b)}{\Gamma(c-a)\Gamma(c-b)}$ at $z=1$ \cite{RakhaElSedy}.  The series converges conditionally for $-1<{\rm Re}(c-a-b)\leq0$ (except at $z=1$), and it does not converge for ${\rm Re}(c-a-b)\leq-1$ \cite{AbramowitzStegun}.

The Gauss hypergeometric function satisfies the differential equation \eqref{2f1de} whenever none of $c$, $c-a-b$, or $a-b$ is an integer \cite{AbramowitzStegun}.  When $c\in\mathbb{Z}^{-}\cup\{0\}$, one uses the fact that the function $\mathbf{F}(a,b;c;z)=\sum_{j=0}^{\infty}\frac{(a)_{j}(b)_{j}}{\Gamma(c+j)}\frac{z^{j}}{j!}$ is a solution for $\left|z\right|<1$ and any $c\in\mathbb{C}$. Most methods that we describe can be modified to compute $_{2}F_{1}(a,b;c;z)$ as well.  The differential equation \eqref{2f1de} has three singular points: $z=0$, $z=1$, and $z=\infty$.

As discussed in Ref.~\cite{Olver1}, one constructs a branch cut between $z=1$ and $z=+\infty$. The branch in the sector $\left|\arg (1-z)\right|<\pi$ is considered to be the \textbf{principle branch}, and we aim to compute $\mathbf{F}$ in this branch.

As was the case with the confluent hypergeometric function, the Gauss hypergeometric function arises in a wide range of applications. For instance, $\mathbf{F}$ can be used to describe transonic adiabatic flow over a smooth bump in an ideal compressible fluid \cite{ChiocchiaGabutti}.  It has also appeared in investigations of a plasma dispersion function \cite{MaceHellberg}, penetration by electrons of a potential barrier \cite{Eckart}, the density of infected nodes in the susceptible-infective-suseptible (SIS) epidemic model on networks \cite{PastorSatorrasVespignani}, scintillation indices of beams \cite{Efimov,Korotkova}, stochastic dynamical systems \cite{WDLL}, electroosmosis of non-Newtonian fluids \cite{ZhaoYang}, and more.


\subsection{Taylor Series} \label{sec:2F1_taylor}

In this section, we carry out computations of the Gauss hypergeometric function using its Taylor series representation \eqref{2F1} and discuss the accuracy of these computations for different parameter and variable regimes.

As with the confluent hypergeometric function in Section \ref{sec:1F1_taylor}, we use two methods to compute the series
\begin{align} \label{2F1seriesdef} 
	_{2}F_{1}(a,b;c;z)=\sum_{j=0}^{\infty}\underbrace{\frac{(a)_{j}(b)_{j}}{(c)_{j}}~\frac{1}{j!}z^{j}}_{C_{j}}\,.
\end{align}

\vspace{0.75em}
\underline{\textbf{Method (a):}} ~We compute

\begin{align*}
	C_{0}&=1\,,\quad S_{0}=C_{0}\,, \\
	C_{j+1}&=C_{j}~\frac{(a+j)(b+j)}{c+j}~\frac{z}{j+1}\,,\quad S_{j+1}=S_{j}+C_{j+1}\,,\quad{}j=0,1,2,\ldots~\,,
\end{align*}
where $C_{j}$ denotes the $(j+1)$th term of the Taylor series \eqref{2F1} and $S_{j}$ denotes the sum of the first $j+1$ terms.

We stop the summation either when $\frac{\left|C_{N+1}\right|}{\left|S_{N}\right|}<\text{\url{tol}}$ for some $\text{\url{tol}}$ and some $N$ or when two successive terms are small compared to $S_{N}$. We then return $S_{N}$ as the solution.

\vspace{0.75em}
\underline{\textbf{Method (b):}} As with the recommended method of Ref.~\cite{Muller} for computing $_{1}F_{1}$ that we discussed in Section \ref{sec:1F1_taylor}, one can compute a recurrence relation for an approximation of the Gauss hypergeometric function in terms of the two previous approximations in an iterative scheme:
\begin{align*}
	\ S_{-1}&{}=S_{0}=1\,,\quad\quad\quad\quad\quad\quad\quad~~{}S_{1}=1+\frac{ab}{c}z\,, \\
	\ r_{j}&{}=\frac{(a+j-1)(b+j-1)}{j(c+j-1)},\quad{}S_{j}=S_{j-1}+(S_{j-1}-S_{j-2})r_{j}z,\quad{}j=2,3,\ldots~\,.
\end{align*}
We stop the summation either when $\frac{\left|S_{N+1}-S_{N}\right|}{\left|S_{N}\right|}<\text{\url{tol}}$ for some $\text{\url{tol}}$ and some $N$ or when two or more successive terms are small.  We then return $S_{N}$ as the solution.

Methods (a) and (b) both amount to truncating the series
\begin{align}
	 S_{\infty}=\sum_{j=0}^{\infty}\frac{(a)_{j}(b)_{j}}{(c)_{j}}~\frac{z^{j}}{j!}\,.
\end{align}
\vspace{0.75em}

By using one of these above methods to compute $F(a,b;c;z)$, we can subsequently compute $\mathbf{F}(a,b;c;z)$ using the expression $F(a,b;c;z)=\Gamma(c)\mathbf{F}(a,b;c;z)$. The exception to this is when $c\in\mathbb{Z}^{-}\cup\{0\}$, in which case more care must be taken (as was also the case for the confluent hypergeometric function).

\begin{figure}[h]
\centering
\includegraphics[width=8.5cm]{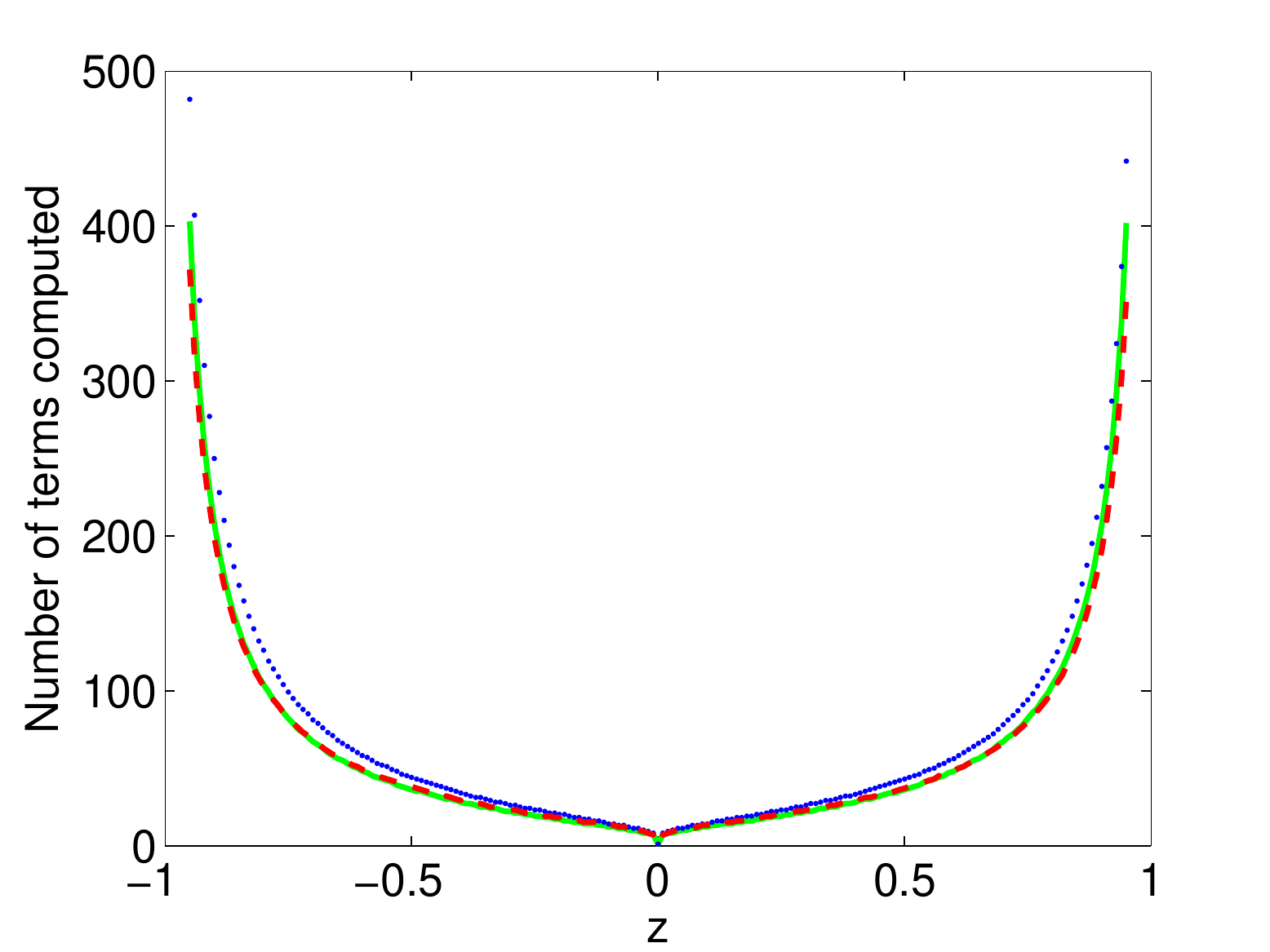}
\caption{(Color online) Number of terms that need to be computed using Taylor series method (a) for evaluating $_{2}F_{1}$ (and hence $\mathbf{F}$) for real $z\in[-0.95,0.95]$ when $a=1.5$, $b=1+2i$, $c=4.5+5i$ (red, dashed), $a=0.15$, $b=0.2$, $c=1.1$ (green, solid), and $a=3$, $b=2$, $c=6.5$ (blue, dotted) before the stopping criterion (two successive small terms) is satisfied.  The method gives results to 14-digit accuracy.
}
\label{fig:figure4.2}
\end{figure}

As shown in Table \ref{appbtable2}, methods (a) and (b) have similar levels of effectiveness in terms of accuracy and number of terms required for computation. Both methods work very successfully when the parameter values are of small magnitude (especially when $|a|, |b|, |c|\leq20$), and the time it takes for the computation is significantly shorter than that taken by the inbuilt {\scshape Matlab} program.  However, computing $\mathbf{F}$ (before the stopping criterion is met) with Taylor series requires more terms than was the case for computing $\mathbf{M}$.  In Fig.~\ref{fig:figure4.2}, we show the number of points required for the computation of three test cases with real $z\in[-1,1]$. As one can see, this illustrates that many more points are required for the computation as $z$ approaches the unit disk.

We can compute $_{2}F_{1}$ very accurately for a range of cases with large real parameter values (with $c<a,b<0$). However, the Taylor series method is much less effective in cases in which either $\left|a\right|$ or $\left|b\right|$ is much greater than $\left|c\right|$.  For such cases, one should employ the recurrence relations that we will discuss in Section \ref{sec:2F1_largeparams_recurrences}.  We conclude that the region in which the Taylor series methods seem to be effective is $\left|z\right|\lesssim 0.9$, provided the values of $a$, $b$, and $c$ do not result in cancellation.

Other series expressions for $_{2}F_{1}$ are also available. For example, Ref.~\cite{Potts} discussed an expansion in terms of $\widehat{z}:=\frac{\sqrt{1-z}-1}{\sqrt{1-z}+1}$.  This expansion converges for $|\widehat{z}|<1$, which is a less restrictive condition than the convergence criterion for \eqref{2F1seriesdef}. As for $\mathbf{M}$, one can exploit sequence acceleration methods for computing the series \eqref{2F1seriesdef}.  For example, one can use the NAG Library routine \url{c06ba} \cite{NAGLibrary}.


\subsubsection{Writing the Gauss Hypergeometric Function as a Single Fraction} \label{sec:2F1_singlefraction}

As for the confluent hypergeometric function, we can compute the Gauss hypergeometric function by representing it as a single fraction. We now examine the accuracy and robustness of this approach and indicate parameter regimes for which it is particularly effective.

\vspace{0.75em}
\underline{\textbf{Method (c):}} As we discussed for the confluent hypergeometric function in Section \ref{sec:1F1_singlefraction}, one can express $_{2}F_{1}(a,b;c;z)$ as a single fraction using recurrence relations.  One starts with $\alpha_{0}=0$, $\beta_{0}=1$, $\gamma_{0}=1$, and $S_{0}=1$, and the recurrence relation is given by
\begin{align}
	\notag \alpha_{j}&=(\alpha_{j-1}+\beta_{j-1})j(c+j-1)\,, \\
	 \notag \beta_{j}&=\beta_{j-1}(a+j-1)(b+j-1)z\,, \\
	 \notag \gamma_{j}&=\gamma_{j-1}j(c+j-1)\,, \\
 \label{nardineq15} S_{j}&=\frac{\alpha_{j}+\beta_{j}}{\gamma_{j}}\,, \quad j=1,2,\ldots\,,
\end{align}
This generates a sequence of approximations $S_{j}$ ($j=1,2, \ldots)$ to $_{2}F_{1}(a,b;c;z)$. We can select stopping criteria in the same way as for the other Taylor series methods discussed above. 
\vspace{0.75em}

As indicated in Table \ref{appbtable2}, this approach works well for small values of the parameters and variable (as a rough guide, it is good for $\left|a\right|,\left|b\right|,\left|c\right|\lesssim20$ and $\left|z\right|\lesssim0.9$). In particular,  this method becomes more successful for computing $\mathbf{F}$ as $\left|c\right|$ gets smaller or $c$ gets closer to an integer.  The reason is the same as it was for the analogous method for confluent hypergeometric functions when computing $\mathbf{M}$ for small $\left|b\right|$ or $b$ close to $-m$, where $m\in\mathbb{Z}^{+}\cup\{0\}$ (see Section \ref{sec:1F1_singlefraction}).

However, this method does struggle more than the Taylor series methods when either $a$ or $b$ has large magnitude (roughly greater than 50), due to a greater risk of overflow---i.e., the program is attempting to compute values that are larger than it can handle---due to the potentially large numerators and denominators in \eqref{nardineq15}. For such cases, one should use other methods, such as the recurrence relations that we will describe in Section \ref{sec:2F1_largeparams_recurrences}.


\subsection{Quadrature Methods} \label{sec:2F1_quadrature}

In Section \ref{sec:1F1_quadrature}, we discussed quadrature methods for computing $\mathbf{M}(a;b;z)$.  We now apply the method of Gauss-Jacobi quadrature to compute $\mathbf{F}(a,b;c;z)$ when ${\rm Re}(c)>{\rm Re}(b)>0$ and $\left|\arg(1-z)\right|<\pi$.  In this parameter and variable regime, the function $\mathbf{F}(a,b;c;z)$ has the integral representation \cite{AbramowitzStegun}
\begin{align} \label{2F1int} 
	\mathbf{F}(a,b;c;z) = \frac{1}{\Gamma(b)\Gamma(c-b)}\int_{0}^{1}(1-zt)^{-a}(1-t)^{c-b-1}t^{b-1}{\rm d}t\,.
\end{align}
Note that $a$ and $b$ can be interchanged in the series definition \eqref{2F1} of $_{2}F_{1}$, so we can also apply equation (\ref{2F1int}) to $\mathbf{F}(b,a;c;z)$ if ${\rm Re}(c)>{\rm Re}(a)>0$.

Transforming $t\mapsto\frac{1}{2}\widetilde{t}+\frac{1}{2}$ and defining Jacobi parameters $\widetilde{\alpha}=c-b-1$ and $\widetilde{\beta}=b-1$ yields \cite{Gautschi3}
\begin{align*}
	\mathbf{F}(a,b;c;z) &= \frac{1}{2^{c-1}\Gamma(b)\Gamma(c-b)}\int_{-1}^{1}\left(\Big(1-\frac{z}{2}\Big)-\frac{1}{2}z\widetilde{t}\right)^{-a}(1-\widetilde{t})^{c-b-1}(1+\widetilde{t})^{b-1}{\rm d}\widetilde{t} \\
\ &=\frac{1}{\Gamma(b)\Gamma(c-b)}\left[\sum_{j=1}^{N_{\mbox{\scriptsize mesh}}}w_{j}^{GJ}\left(\Big(1-\frac{z}{2}\Big)-\frac{1}{2}zt_{j}^{GJ}\right)^{-a}+E_{N_{\mbox{\scriptsize mesh}}}(a;b;z)\right]\,,
\end{align*}
where $t_{j}^{GJ}$ and $w_{j}^{GJ}$ are the Gauss-Jacobi nodes and weights on $[-1,1]$ and $N_{\mbox{\scriptsize mesh}}$ is the number of mesh points. Error bounds for this method are discussed in Ref.~\cite{Gautschi3}.

If ${\rm Re}(c)>{\rm Re}(a)>0$, then switching the parameters $a$ and $b$ in the definition of $\mathbf{F}(a,b;c;z)$ allows one to apply the method of Gauss-Jacobi quadrature. As was the case for $\mathbf{M}$, the integrand in \eqref{2F1int} blows up at the end points of the integral. This motivates the choice of Gauss-Jacobi quadrature to perform the required integration numerically.

The results in Table \ref{appbtable2} illustrate that applying Gauss-Jacobi quadrature to the integral in \eqref{2F1int} is a useful method for computing the Gauss hypergeometric function when ${\rm Re}(c)>{\rm Re}(b)>0$. If the moduli of the parameters $a$, $b$, and $c$ are at least 50--100, then it is an improvement to the Taylor series approaches. The method also works well near $z = e^{\pm i\pi/3}$, which are notoriously difficult to compute (see the discussions in Sections \ref{sec:2F1_unitdisc} and \ref{sec:2F1_analyticcont}), especially for smaller values of $|a|$ and $|b|$.

Therefore, as with $_{1}F_{1}$ (see Section \ref{sec:1F1_quadrature}), applying Gauss-Jacobi quadrature to compute $_{2}F_{1}$ is potentially a useful method when the parameters $a$, $b$, and $c$ have moduli that are not too large. However, due to the relatively high computational cost of generating the quadrature nodes and weights, we only recommend this approach be applied in parameter regimes where a series method is likely to generate inaccurate solutions.


\subsection{Computing $\mathbf{F}$ for $z$ Near or Outside of the Unit Disk} \label{sec:2F1_unitdisc}

The methods that we have discussed thus far are generally very useful for computing $\mathbf{F}$ for $|z|\lesssim 0.9$.  In this section, we discuss how to use these computations to compute $\mathbf{F}$ effectively for any $z\in\mathbb{C}$.  This can be done either by applying transformation formulas to the results from previous sections or by using expansions derived for the special cases $b-a\in\mathbb{Z}$ and $c-a-b\in\mathbb{Z}$.


\subsubsection{Transformation Formulas} \label{sec:2F1_unitdisc_transformations}

Because the series \eqref{2F1} converges only for $\left|z\right|<1$, and because it converges more rapidly as $\left|z\right|$ becomes smaller, it is important to use transformation formulas that reduce the problem of computing \eqref{2F1} for a value of $\left|z\right|$ near or greater than $1$ to a problem of computing the series for a new variable $w$ with much smaller magnitude. We describe such transformation formulas in this section.

The idea of these transformations is to map as large a region of the complex plane as possible onto disks $\left|w\right|\leq\rho$ for a 
number $\rho \in (0,1]$ that is preferably as close to $0$ as possible. This is desirable because the function $\mathbf{F}$ can be computed faster and more accurately when $\left|z\right|$ is near $0$.  Finding representations of $\mathbf{F}$ that make it possible to carry out the computation in terms of the new variable $w$ allows one to obtain more accurate results than by using the methods that we described previously. We show such transformation formulas for real $z$ in Table \ref{transtable}.  These transformations map any $z\in\mathbb{R}$ to a new variable $w\in[0,\frac{1}{2}]$.

\begin{table}[h]
\renewcommand{\arraystretch}{1.15}
\begin{center}
\begin{tabular}{|c|c|c|}
\hline
Case & Interval & Transformation \\ \hline \hline
1 & $-\infty<z<-1$ & $w=\frac{1}{1-z}$ \\
2 & $-1\leq z<0$ & $w=\frac{z}{z-1}$ \\
3 & $0\leq z\leq\frac{1}{2}$ & $w=z$ \\
4 & $\frac{1}{2}<z\leq1$ & $w=1-z$ \\
5 & $1<z\leq2$ & $w=1-\frac{1}{z}$ \\
6 & $2<z<\infty$ & $w=\frac{1}{z}$ \\
\hline
\end{tabular}
\caption{List of transformations of $z\in\mathbb{R}$ \cite{Forrey} for which $0\leq w\leq\rho=\frac{1}{2}$.}
\label{transtable}
\end{center}
\end{table}

\begin{figure}
  \centering
  \subfloat{\label{fig:rho0.6}\includegraphics[width=0.5\textwidth]{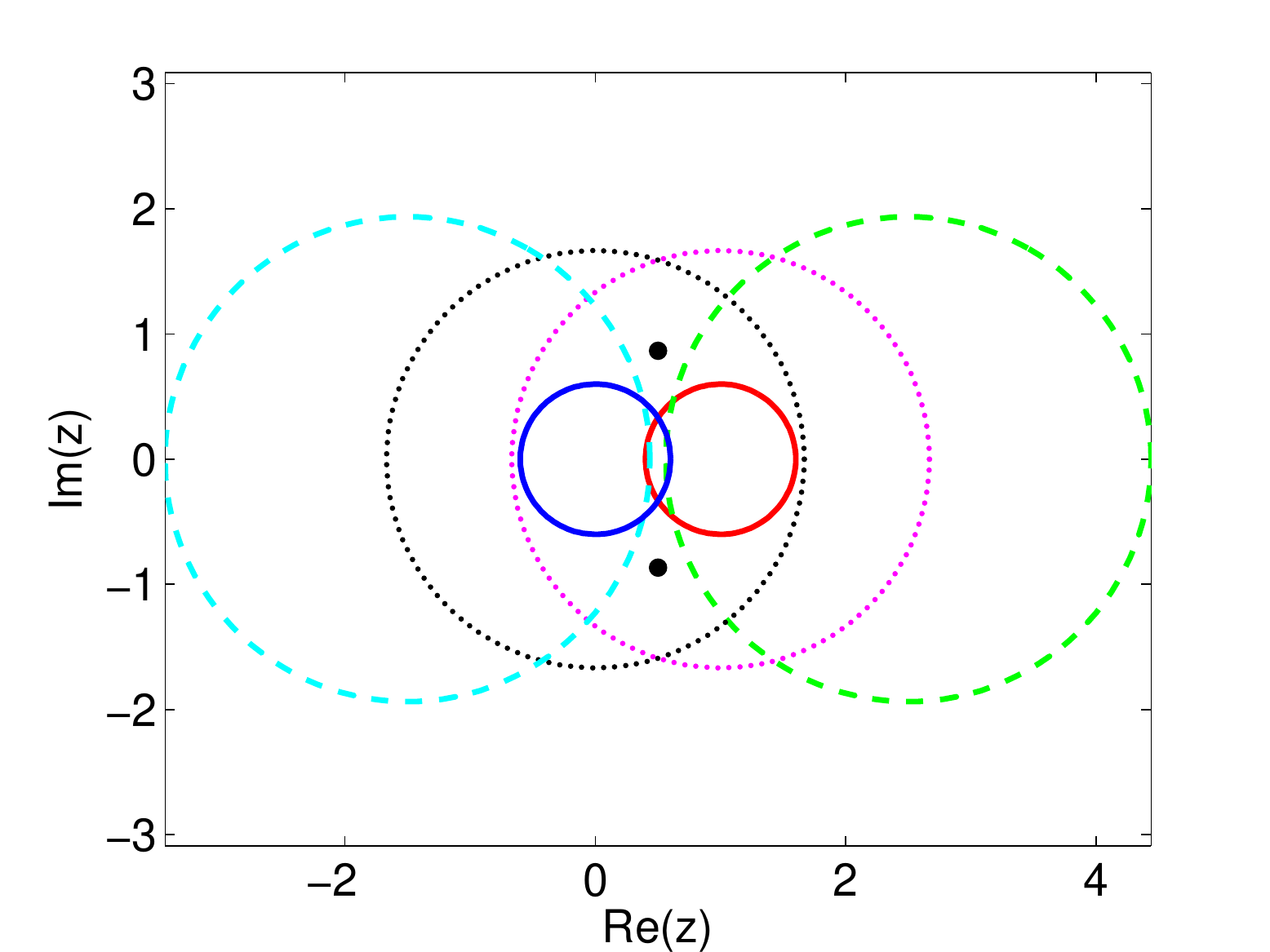}}
  \subfloat{\label{fig:rho0.8}\includegraphics[width=0.5\textwidth]{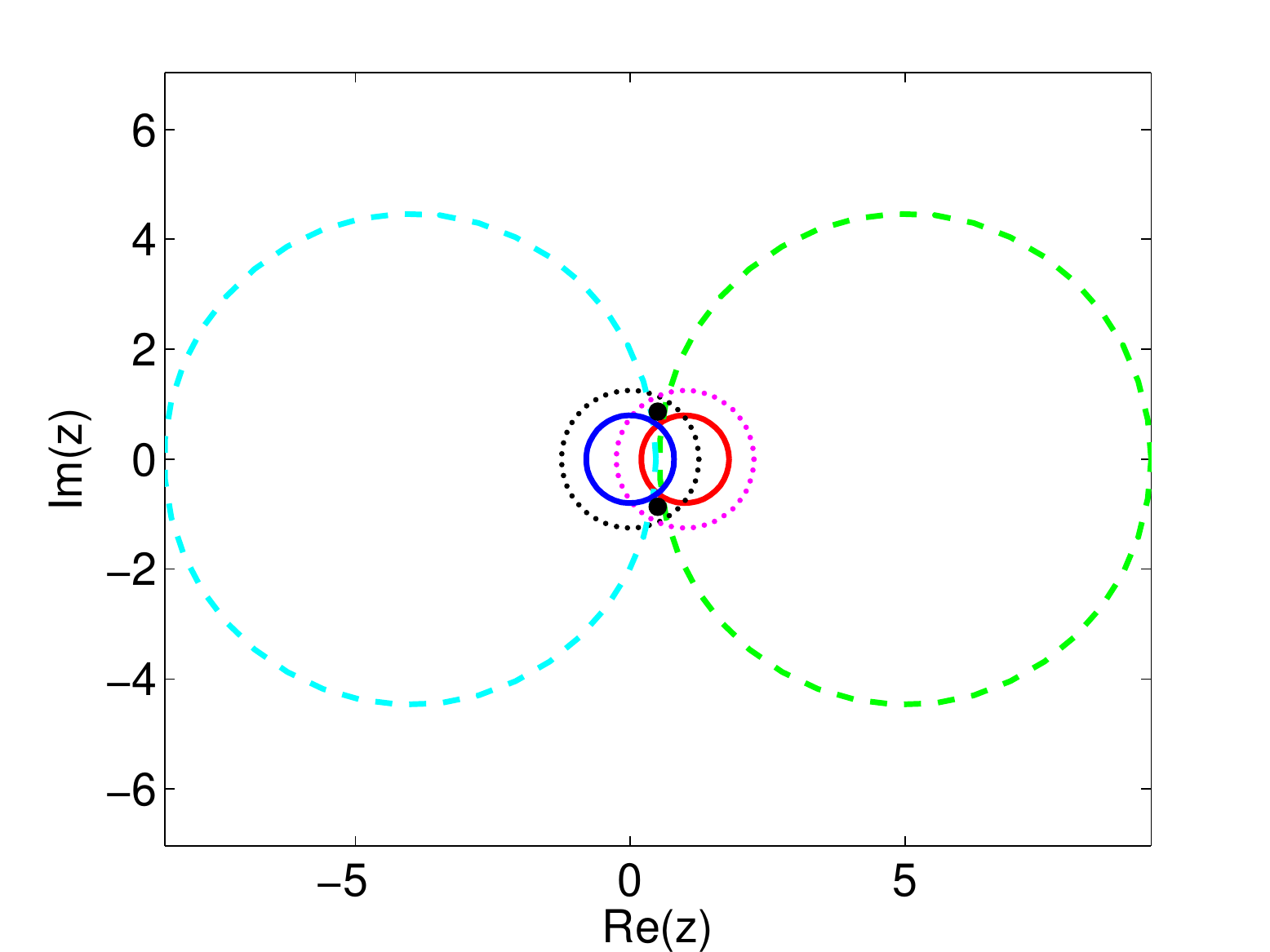}}
  \caption{(Color online) The curves $\left|z\right|=\rho$ (dark blue, solid, left inner), $\left|\frac{1}{z}\right|=\rho$ (black, dotted, left middle), $\left|1-z\right|=\rho$ (red, solid, right inner), $\left|\frac{1}{1-z}\right|=\rho$ (purple, dotted, right middle), $\left|\frac{z}{z-1}\right|=\rho$ (sky blue, dashed, left outer), and $\left|1-\frac{1}{z}\right|=\rho$ (green, dashed, right outer), for $\rho=0.6$ (left panel) and $\rho=0.8$ (right panel). We show the points $z=e^{\pm{}i\pi/3}$ using black dots. [These figures are adapted from plots in Ref.~\cite{GST1,GST3}.]
  }
  \label{fig:graphs}
\end{figure}

For complex $z$, the problem is more complicated. In Fig.~\ref{fig:graphs}, we show plots of $\left|w\right|=\rho$ for each of the six expressions for $w$ given in Table \ref{transtable} for $\rho=0.6$ and $\rho=0.8$. If one wishes to apply one of the six transformations in Table \ref{transtable}, it is required that $\left|w\right|\leq\rho$ be satisfied for at least one representation of $w$ in the table.  The region in Fig.~\ref{fig:graphs} in which none of the representations of $w$ satisfy $\left|w\right|\leq\rho$ either surround or are near the points $z=e^{\pm i\pi/3}=\frac{1}{2}(1\pm i\sqrt{3})$ (marked as dots). As $\rho$ increases towards the value $1$, the region in which none of the transformations satisfy $\left|w\right|\leq\rho$ gets smaller, but it remains near the points $z=e^{\pm i\pi/3}$ due to the fact that the set $\{e^{i\pi/3},e^{-i\pi/3}\}$ is mapped to itself by each of the six transformations in Table \ref{transtable}.  We discuss the case $z\approx e^{\pm i\pi/3}$ in more detail in 
Section \ref{sec:2F1_analyticcont} and noted it briefly in Section \ref{sec:2F1_quadrature}. The five transformations that are useful when $z \not\approx e^{\pm{}i\pi/3}$ are \cite{AbramowitzStegun,BeckenSchmelcher,EMOT,Forrey,ZhangJin}
\begin{align} \label{treq1} 
	\frac{\sin(\pi[b-a])}{\pi}~\mathbf{F}(a,b;c;z) &{}={}\frac{(1-z)^{-a}}{\Gamma(b)\Gamma(c-a)}~
\mathbf{F}\left(a;c-b;a-b+1;\frac{1}{1-z}\right) \\
\ \nonumber &+\frac{(1-z)^{-b}}{\Gamma(a)\Gamma(c-b)}~
\mathbf{F}\left(b;c-a;b-a+1;\frac{1}{1-z}\right)\,, \\
	 \label{treq2} 
\mathbf{F}(a,b;c;z) &{}={}(1-z)^{-a}~
\mathbf{F}\left(a;c-b;c;\frac{z}{z-1}\right)\,, \\
	 \label{treq3} 
\frac{\sin(\pi[c-a-b])}{\pi}~\mathbf{F}(a,b;c;z) &{}={}\frac{1}{\Gamma(c-a)\Gamma(c-b)}~
\mathbf{F}(a;b;a+b-c+1;1-z) \\
 \nonumber &+\frac{(1-z)^{c-a-b}}{\Gamma(a)\Gamma(b)}~
\mathbf{F}(c-a;c-b;c-a-b+1;1-z)\,, \\
	 \label{treq4} 
\frac{\sin(\pi[c-a-b])}{\pi}~\mathbf{F}(a,b;c;z) & {}={}\frac{z^{-a}}{\Gamma(c-a)\Gamma(c-b)}~
\mathbf{F}\left(a;a-c+1;a+b-c+1;1 - \frac{1}{z}\right) \\
 \nonumber &+\frac{z^{a-c}\left(1-z\right)^{c-a-b}}{\Gamma(a)\Gamma(b)}~\mathbf{F}\left(c-a;1-a;c-a-b+1;1 - \frac{1}{z}\right)\,, \\
 	\label{treq5} 
\frac{\sin(\pi[b-a])}{\pi}~\mathbf{F}(a,b;c;z) &{}={}\frac{(-z)^{-a}}{\Gamma(b)\Gamma(c-a)}~
\mathbf{F}\left(a;a-c+1;a-b+1;\frac{1}{z}\right) \\
 \nonumber &+\frac{(-z)^{-b}}{\Gamma(a)\Gamma(c-b)}~
\mathbf{F}\left(b-c+1;b;b-a+1;\frac{1}{z}\right)\,,
\end{align}
where \eqref{treq1} and \eqref{treq3}--\eqref{treq5} require that $\left|\arg(1-z)\right|<\pi$ and \eqref{treq4}--\eqref{treq5} additionally require that $\left|\arg{}z\right|<\pi$.

We tested these formulas on a large range of parameter and variable values and found that they successfully computed $\mathbf{F}$ for variable values close to the unit disk or with modulus greater than $1$ in terms of new variables with smaller magnitude, provided $z$ was not too close to $e^{\pm{}i\pi/3}$. However, due to the presence of $\Gamma(a-b)$, $\Gamma(b-a)$, $\Gamma(c-a-b)$, and $\Gamma(a+b-c)$ terms in the numerators of \eqref{treq1}--\eqref{treq5}, the cases $b-a\in\mathbb{Z}$ and $c-a-b\in\mathbb{Z}$ cannot be handled by applying the formulas in Table \ref{transtable}. There are also numerical issues when either $b-a$ or $c-a-b$ is close to an integer, as two large conflicting terms are being added together (leading to cancellation).


\subsubsection{Cases $b-a\in\mathbb{Z}$ and $c-a-b\in\mathbb{Z}$} \label{sec:2F1_unitdisc_integers}

References \cite{AbramowitzStegun,EMOT} discussed formulas that avoid this cancellation issue when either $b-a$ or $c-a-b$ is exactly equal to an integer. These cases can be computed using the ideas detailed in Section \ref{sec:1F1_taylor}.  Reference \cite{Forrey} derived expressions for $b-a$ or $c-a-b$ close to an integer using a polynomial expansion for the Gamma function.  Reference \cite{MichelStoitsov} showed how to exploit the Lanczos expansion and other properties of the Gamma function for successful computation when $b-a$ or $c-a-b$ is close to an integer.


\subsection{Analytic Continuation Formulas for $z \approx e^{\pm i\pi/3}$} \label{sec:2F1_analyticcont}

Applying the transformation formulas of Section \ref{sec:2F1_unitdisc} becomes a much less viable method as $z$ approaches the points $z=e^{\pm i\pi/3}=\frac{1}{2}(1\pm i\sqrt{3})$. The reason is as follows: whatever value $\rho \in (0,1]$ is taken for $\left|z\right|\leq\rho$, it is not possible to map the points $z=e^{\pm i\pi/3}$ onto $w$ within a disk of radius less than $1$, as points are mapped to themselves or each other under any of the six transformations in Table \ref{transtable}.  As $\rho$ increases towards $1$, an increasingly large number of points close to $z=e^{\pm i\pi/3}$ are mapped into such a disk, but the points themselves can never be mapped into a disk, and points very close to them require a computation that involves a prohibitively large value of $\rho$ (for which the methods discussed thus far do not generate an accurate result).

To resolve this issue, Refs.~\cite{Buhring1,Buhring2} gave an expansion in the form of a \emph{continuation formula}:
\begin{align} \label{continuation1} 
	\mathbf{F}(a,b;c;z)={}&\frac{\Gamma(b-a)}{\Gamma(b)\Gamma(c-a)}(z_{0}-z)^{-a}\sum_{j=0}^{\infty}d_{j}(a,z_{0})(z-z_{0})^{-n} \\
\ \notag &+{}\frac{\Gamma(a-b)}{\Gamma(a)\Gamma(c-b)}(z_{0}-z)^{-b}\sum_{j=0}^{\infty}d_{j}(b,z_{0})(z-z_{0})^{-n}\,,
\end{align}
where $\left|\arg(z_{0}-z)\right|<\pi$, the difference $b-a$ is not an integer, and $d_{j}$ is defined by
\begin{align}
	 d_{-1}(\upsilon,z_{0})={}&0\,,~d_{0}(\upsilon,z_{0})=1\,, \\
	 \notag d_{j}(\upsilon,z_{0})={}&\frac{j+\upsilon-1}{j(j+2\upsilon-a-b)}\left[\{(j+\upsilon)(1-2z_{0})+(a+b+1)z_{0}-c\}d_{j-1}(\upsilon,z_{0})\right. \\
\ \notag &~~~~~~~~~~~~~~~~~~~~~~~~\left.+z_{0}(1-z_{0})(j+\upsilon-2)d_{j-2}(\upsilon,z_{0})\right]\,,~~j=1,2,\ldots~\,.
\end{align}

The series in \eqref{continuation1} converges at every point outside of the disk $\left|z-z_{0}\right|=\max\{\left|z_{0}\right|,\left|z_{0}-1\right|\}$, so one must choose an appropriate $z_{0}$ to implement this method. The case noted in Ref.~\cite{Buhring1} is $z_{0}=\frac{1}{2}$, which results in convergence outside of the disk $\left|z-\frac{1}{2}\right|=\frac{1}{2}$ and, in particular, includes points close to $e^{\pm i\pi/3}$. 
For further details about this method, see Refs.~\cite{Buhring1,Buhring2,Kalla}.

For $b-a\in\mathbb{Z}$, the expansion \eqref{continuation1} is not valid, as one of the Gamma functions in the numerators of the two terms in \eqref{continuation1} is infinitely large. Reference \cite{Buhring1} discussed these issues and includes an alternative expression for $b-a=0$.  Reference \cite{Luke4} used a limiting process for $b-a\in\mathbb{Z}\backslash\{0\}$.

The results in Table \ref{appbtable2} suggest that using the expansion in \eqref{continuation1} is very accurate in a region that troubled the methods from Sections \ref{sec:2F1_taylor}--\ref{sec:2F1_unitdisc}. Using the continuation formula \eqref{continuation1} is effective not only when $z$ is equal or close to $e^{\pm i\pi/3}$ but also when $z$ lies outside of the unit disk (see the fourth and fifth rows of Table \ref{appbtable2}). However, it is ineffective when $\left|a\right| \gtrsim 30$, $\left|b\right|\gtrsim30$, or $\left|c\right|\gtrsim70$, but this problem can be resolved by applying recurrence relations from Section \ref{sec:2F1_largeparams_recurrences}.

Reference \cite{MichelStoitsov} presents an alternative method to that suggested in Ref.~\cite{Buhring1} for computing $_{2}F_{1}(a,b;c;z)$ near the points $z=e^{\pm{}i\pi/3}$.  The idea is again to expand about a point $z_{0}$ at which the expansion can be computed more easily, but this time one uses a Taylor expansion about this point in conjunction with the recurrence relation corresponding to the hypergeometric differential equation \eqref{2f1de}. The expansion that was obtained in Ref.~\cite{MichelStoitsov} is
\begin{align*}
	_{2}F_{1}(a,b;c;z) = \sum_{n=0}^{\infty}q_{n}(z-z_{0})^{n}\,,
\end{align*}
where
\begin{align*}
	q_{0}&{}={} _{2}F_{1}(a,b;c;z_0)\,,\quad{}q_{1}=\frac{\rm d}{{\rm d}z_{0}}\left[
_{2}F_{1}(a,b;c;z_0)\right]=\frac{ab}{c}~_{2}F_{1}(a+1;b+1;c+1;z_0)\,, \\
	 q_{n}&{}=\frac{1}{z_{0}(1-z_{0})(n+2)}\left[\Big(n(2z_{0}-1)-c+(a+b+1)z_{0}\Big)q_{n+1}+\frac{(a+n)(b+n)}{n+1}q_{n}\right]\,,
\end{align*}
for $n=0,1,2, \ldots$. One can also use these relations to compute the function $\mathbf{F}(a,b;c;z)$. The quantities $q_{0}$ and $q_{1}$ should be calculated using a method appropriate for the parameter regimes in question.  As in Ref.~\cite{MichelStoitsov}, we take
\begin{align*}
	z_{0}=r_{0}\frac{z}{|z|}\,,\quad{}r_{0}=\left\{\begin{array}{rl}
0.9 & \text{if}~|z|\leq 1\,,\\
1.1 & \text{if}~|z|>1\,.\\
\end{array}~~~~\right.
\end{align*}

The method is advantageous compared to that in Ref.~\cite{Buhring1} because there are no Gamma functions in the expression, so no singularities are encountered when $b-a\in\mathbb{Z}$ or $c-a-b\in\mathbb{Z}$.  As discussed in Ref.~\cite{MichelStoitsov}, this method gives better accuracy for large $|a|$ and $|b|$.  However, the method can be slow computationally due to the work encountered when computing $q_{0}$ and $q_{1}$.


\subsection{Techniques for Large Values of Parameters $|a|$, $|b|$, and $|c|$} \label{sec:2F1_largeparams}

In this section, we detail two widely researched techniques that can be useful for evaluating $\mathbf{F}$ when $|a|$, $|b|$, and $|c|$ are sufficiently large that the function cannot be evaluated effectively using methods that we described previously. The two methods we now discuss entail the use of recurrence relations---involving computing $\mathbf{F}$ for ``nicer'' parameter values and applying the recurrence relations to obtain a solution for more awkward values---or uniform asymptotics, which involves computing series for $\mathbf{F}$ in terms of other special functions.


\subsubsection{Recurrence Relations} \label{sec:2F1_largeparams_recurrences}

As we discussed for the confluent hypergeometric function in Section \ref{sec:1F1_recurrences}, we aim to overcome the lack of accuracy from attempting to compute the Gauss hypergeometric function for large values of $\left|{\rm Re}(a)\right|$, $\left|{\rm Re}(b)\right|$, and $\left|{\rm Re}(c)\right|$---or, as mentioned in Section \ref{sec:2F1_summary}---which can sometimes occur when $\left|{\rm Re}(c)\right|$ is too small. We explain how this problem can be addressed using recurrence relations. We then reformulate the computation as one involving values of $\left|{\rm Re}(a)\right|$, $\left|{\rm Re}(b)\right|$, and $\left|{\rm Re}(c)\right|$ closer to zero, which can then be computed accurately by the methods discussed previously.

Recurrence relations involving the Gauss hypergeometric function $\mathbf{F}(a,b;c;z)$ are discussed in Refs.~\cite{FLS,GST1,GST2,IbrahimRakha,Temme1}. The four recurrence relations that we discuss in the present paper are described in Ref.~\cite{GST1}.  Following the usual notation in the literature, we denote them by $(++0)$, $(00+)$, $(++-)$, and $(+0-)$.  We state them in Appendix \ref{sec:appendixd}.

As with Section \ref{sec:1F1_recurrences} for $\mathbf{M}$, we seek minimal solutions of the recurrence relations.  Unlike for $\mathbf{M}$, however, the four recurrence relations for $\mathbf{F}$ have different minimal solutions in different regions of the complex plane. We give the minimal solutions for each of the four recurrence relations \cite{GST2} in Table \ref{minsolns2f1}.

\begin{table}[h]
\renewcommand{\arraystretch}{1.15}
\begin{center}
\begin{tiny}
\begin{tabular}{|@{}c@{}|@{}l@{}|}
\hline
Relation and valid region of $\mathbb{C}$ & ~Minimal solution \\ \hline \hline
$(++0)$ & ~$\Gamma(1+a-c+n)\Gamma(1+b-c+n)~
\mathbf{F}(a+n;b+n;1+a+b-c+2n;1-z)$ \\
$\mathbb{C}\backslash\{z\leq0\}$ & \\ \hline
$(00+)$ & ~$\Gamma(c+n)~\mathbf{F}(a;b;c+n;z)$ \\
${\rm Re}(z)<\frac{1}{2}$ & \\ \hline
$(00+)$ & ~$\left(\frac{z-1}{z}\right)^{n}\Gamma(c+n)~
\mathbf{F}(1-a;1-b;1-a-b+c+n;1-z)$ \\
${\rm Re}(z)>\frac{1}{2}$ & \\ \hline
$(++-)$ & ~$\left(\frac{z}{(z-1)^{3}}\right)^{n}\frac{\Gamma(b-c+1+2n)\Gamma(a-c+1+2n)}{\Gamma(a+n)\Gamma(b+n)\Gamma(1-c+n)}~
\mathbf{F}(1-a-n;1-b-n;2-c+n;z)$~ \\
Inside curve in Fig. \ref{fig:recgraphs} (left) & \\ \hline
$(++-)$ & ~$\frac{\Gamma(b-c+1+2n)\Gamma(a-c+1+2n)}{\Gamma(1-c+n)}~\mathbf{F}(a+n;b+n;1+a+b-c+3n;1-z)$ \\
Outside curve in Fig. \ref{fig:recgraphs} (left) & \\ \hline
$(+0-)$ & ~$\left(\frac{-z}{(1-z)^{2}}\right)^{n}\frac{\Gamma(b-c+1+n)\Gamma(a-c+1+2n)}{\Gamma(a+n)\Gamma(1-c+n)}~\mathbf{F}(1-a-n;1-b;2-c+n;z)$ \\
~Inside inner curve in Fig. \ref{fig:recgraphs} (right)~ & \\ \hline
$(+0-)$ & ~$\frac{\Gamma(b-c+1+n)\Gamma(a-c+1+2n)}{\Gamma(1-c+n)}~
\mathbf{F}(a+n;b;1+a+b-c+2n;1-z)$ \\
Between curves in Fig. \ref{fig:recgraphs} (right) & \\ \hline
$(+0-)$ & ~$\left(\frac{-z}{(1-z)^{2}}\right)^{n}\frac{\Gamma(1+a-c+2n)}{\Gamma(1-c+n)}~
\mathbf{F}(1-b;-b+c-n;1+a-b+n;\frac{1}{z})$ \\
Outside curves in Fig. \ref{fig:recgraphs} (right) & \\ \hline
\end{tabular}
\caption{Minimal solutions of the four recurrence relations for $\mathbf{F}$ for $z$ in different regions of the complex plane. See the discussion in Ref.~\cite{GST2}.}
\label{minsolns2f1}
\end{tiny}
\end{center}
\end{table}

\begin{figure}
  \centering
  \subfloat{\label{fig:recgraphs1}\includegraphics[width=0.5\textwidth]{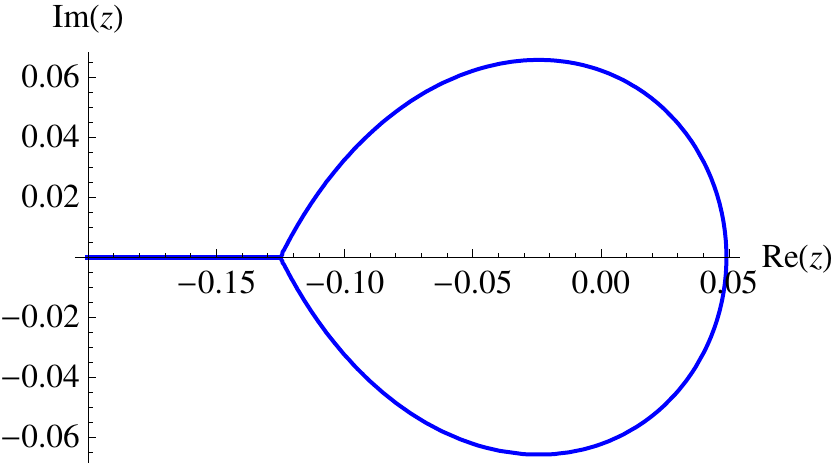}}
  \hspace{2em}
  \subfloat{\label{fig:recgraphs2}\includegraphics[width=0.375\textwidth]{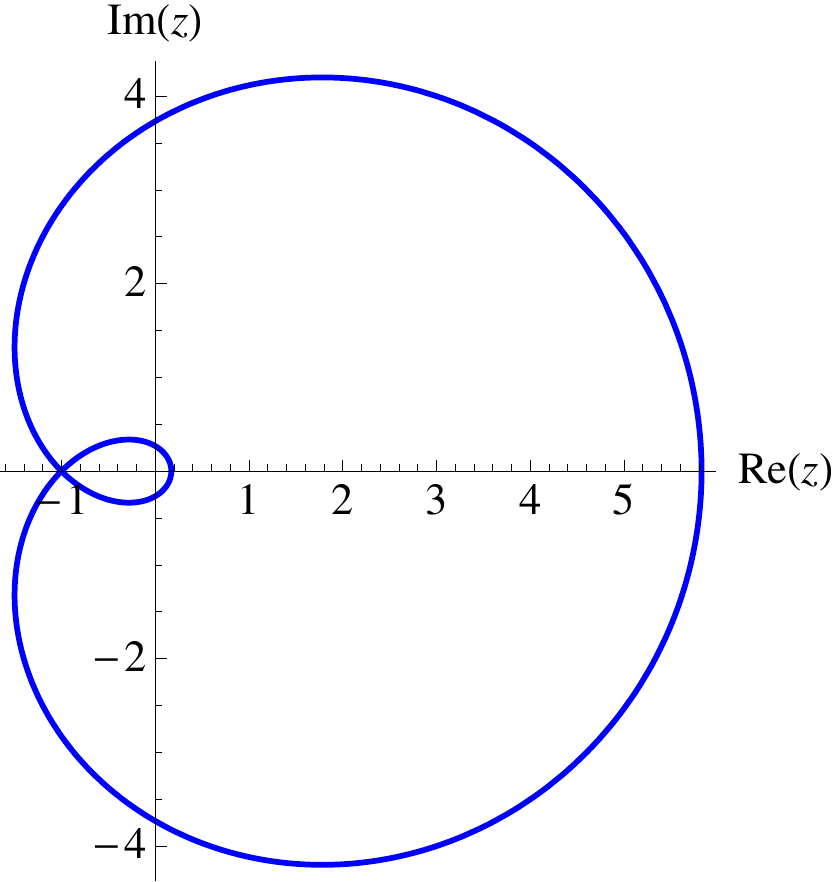}}
  \caption{Relevant regions for minimal solutions of recurrence relations (left) $(++-)$ and (right) $(+0-)$. With reference to Table \ref{minsolns2f1}, on the left plot, ``Inside curve'' refers to the enclosed region and ``Outside curve'' refers to the remainder of $\mathbb{C}$; on the right plot, ``Inside inner curve'' refers to the smaller enclosed region, with ``Between curves'' meaning the larger enclosed region, and ``Outside curves'' the remainder of $\mathbb{C}$.
    }
  \label{fig:recgraphs}
\end{figure}

Therefore, as with $\mathbf{M}$ in Section \ref{sec:1F1_recurrences}, we apply two different methods. First, we take the minimal solutions of the four recurrence relations in specific regions and apply the recurrence relations backwards. Second, we use the minimal solutions of the recurrence relations to apply the recurrence relations forwards using Miller's algorithm or Olver's algorithm.

A viable way of computing a hypergeometric function with parameters whose real parts have large modulus can therefore be to compute two hypergeometric functions whose real parts have smaller modulus, and then apply recurrence relations as appropriate.  However, observe from the minimal solutions of the recurrence relations in Table \ref{minsolns2f1} that this method is constrained by the ability of a programming language to compute the Gamma function for a variable value with large modulus (as was also the case with the recurrence relation techniques of Section \ref{sec:1F1_recurrences} for computing $\mathbf{M}$).  This caveat notwithstanding, the methods discussed in this section are useful for carrying out computations of $\mathbf{F}$ with arguments of large real part.


\subsubsection{Uniform Asymptotic Expansions} \label{sec:2F1_largeparams_uaes}

Another potentially viable method for dealing with large values of the parameters $|a|$, $|b|$, and $|c|$ when computing $\mathbf{F}(a,b;c;z)$ is to exploit uniform asymptotic expansions.  In particular, the following three uniform asymptotic expansions are useful:
\begin{itemize}

\item \underline{\textbf{Expansion for $\mathbf{F}\left(a+\lambda;b-\lambda;c;\frac{1}{2}-\frac{1}{2}z\right)$ as $\lambda\rightarrow\infty$ \cite{Jones}:}}~~For fixed $a,b,c\in\mathbb{C}$ and $|\arg{}z|<\pi$ \cite{Jones} (see also \cite{Dunster1}),
\begin{align*}
_{2}F_{1}\left(a+\lambda;b-\lambda;c;\frac{1}{2}-\frac{1}{2}z\right)\sim{}&2^{\frac{1}{2}(a+b-1)}(z-1)^{-\frac{1}{2}c}(z+1)^{\frac{1}{2}(c-a-b-1)} \\
\ &{}\times\left(\frac{\sinh\zeta}{\zeta}\right)^{1/2}\sum_{m=0}^{\infty}C_{m}(\zeta)\frac{I_{c-1+m}(\alpha\zeta)}{\alpha^{m+c-1}}\,,
\end{align*}
as $\lambda\rightarrow\infty$ for $|\arg{}\lambda|\leq\pi-\delta<\pi$, where $I_{\nu}(z)=i^{-\nu}J_{\nu}(iz)$ and $J_{\nu}(z)$ is the Bessel function (see Appendix \ref{sec:appendixc}). Additionally, $z=\cosh\zeta$, $\alpha=\frac{1}{2}(a-b)+\lambda$, and $C_{j}(\zeta)$ is defined by $C_{0}(\zeta)=1$ and
\begin{align*}
	C_{m}''(\zeta)+\{1-2(\nu+m)\}\frac{C_{m}'(\zeta)}{\zeta}+\{m(m+2\nu)\}\frac{C_{m}(\zeta)}{\zeta^{2}}+2C_{m+1}'(\zeta)=\psi(\zeta)C_{m}(\zeta)\,,
\end{align*}
where $m=0,1,2,\ldots$, the parameter $\nu=\alpha-\frac{1}{2}$, the function $\psi(z)=\frac{\Gamma'(z)}{\Gamma(z)}$ is the psigamma function, and
\begin{align*}
	 \psi(\zeta)&{}=\left(\frac{3}{4}+4c_{1}\right)\left(\frac{1}{\sinh^{2}\zeta}-\frac{1}{\zeta^{2}}\right)+\frac{c_{1}-\frac{1}{2}c_{3}}{\cosh^{2}\frac{1}{2}\zeta}\,, \\
\ c_{1}&{}=\frac{1}{4}\left\{(c-1)^{2}-1\right\}\,,\quad{}c_{3}=\frac{1}{2}\{(a+b-c)^{2}-1\}\,.
\end{align*}

\item \underline{\textbf{Expansion for $_{2}F_{1}\left(a;b-\lambda;c+\lambda;-z\right)$ as $\lambda\rightarrow\infty$ \cite{OldeDaalhuis2}:}}~~For fixed $a,b,c\in\mathbb{C}$ and $|\arg{}z|<\pi$ \cite{OldeDaalhuis2},
\begin{align*}
	_{2}F_{1} &{}\left(a;b-\lambda;c+\lambda;-z\right)\sim\frac{2^{\lambda}(1+z)^{\lambda-a}\Gamma(c+\lambda)\Gamma(1-b+\lambda)}{z^{\lambda/2}\Gamma(c-b+2\lambda)\sqrt{2\pi}} \\
\ &{}\times\left(\lambda^{\frac{1}{2}(a-1)}U\left(a-\frac{1}{2},-\alpha\sqrt{\lambda}\right)\sum_{s=0}^{\infty}\frac{\gamma_{0,s}}{\lambda^{s}}+\lambda^{\frac{1}{2}(a-2)}U\left(a-\frac{3}{2},-\alpha\sqrt{\lambda}\right)\sum_{s=0}^{\infty}\frac{\gamma_{1,s}}{\lambda^{s}}\right)\,,
\end{align*}
as $\lambda\rightarrow\infty$ in $|\arg{}\lambda|\leq\frac{\pi}{2}-\delta<\frac{\pi}{2}$, where $U(a,z)$ is the parabolic cylinder function (see Appendix \ref{sec:appendixc}).  Additionally, $\alpha=-\sqrt{2\text{ln}\frac{(z+1)^{2}}{4z}}$, and the coefficients $\gamma_{j,s}$ are given (for $j=0,1$) by
\begin{align*}
	\gamma_{j,s}=\frac{s!}{(2\pi{}i)^{2}}\oint_{\{0,\tau_{c}\}}\oint_{\{u(\tau)\}}\frac{(u-\alpha)^{a-1}u^{1-j}(\frac{1}{2}-\tau)^{b-1}(\frac{1}{2}+\tau)^{a-c}}{\left(\frac{1}{2}u^{2}+\ln(1-4\tau^{2})\right)^{s+1}(\tau-\tau_{c})^{a}}{\rm d}u{\rm d}\tau\,,
\end{align*}
where $\tau_{c}=\frac{z-1}{2(z+1)}$. Here, the $\tau$-contour is a simple loop encircling the values $0$ and $\tau_{c}$, and, for each $\tau$ on this contour, the $u$-contour is a small loop encircling the point $u(\tau):=\sqrt{-2\text{ln}(1-4\tau^{2})}$ \cite{OldeDaalhuis2}.

\item \underline{\textbf{Expansion for $_{2}F_{1}\left(a+\lambda;b+2\lambda;c;-z\right)$ as $\lambda\rightarrow\infty$ \cite{OldeDaalhuis3}:}}~~For fixed $a,b,c\in\mathbb{C}$ and $|\arg{}z|<\pi$,
\begin{align*}
	 _{2}F_{1}&{}\left(a+\lambda;b+2\lambda;c;-z\right)\sim(1+z)^{-\frac{3}{2}\lambda}\frac{\Gamma(c)\Gamma(1-c+b+2\lambda)}{\Gamma(b+2\lambda)} \\
\ &{}\times\left(\left[e^{\left(a-c+\lambda+\frac{1}{3}\right)\pi{}i}\text{Ai}\left((e^{-\pi{}i}\lambda)^{2/3}x\right)+e^{-\left(a-c+\lambda+\frac{1}{3}\right)\pi{}i}\text{Ai}\left((e^{\pi{}i}\lambda)^{2/3}x\right)\right]\chi_{1}\right. \\
\ &{}\left.-\left[e^{\left(a-c+\lambda+\frac{2}{3}\right)\pi{}i}\text{Ai}'\left((e^{-\pi{}i}\lambda)^{2/3}x\right)+e^{-\left(a-c+\lambda+\frac{2}{3}\right)\pi{}i}\text{Ai}'\left((e^{\pi{}i}\lambda)^{2/3}x\right)\right]\chi_{0}\right)\,,
\end{align*}
where
\begin{align*}
	 \chi_{1}=\sum_{s=0}^{\infty}(-1)^{s}\gamma_{1,s}\lambda^{-s-\frac{1}{3}},\quad\chi_{0}=\sum_{s=0}^{\infty}(-1)^{s}\gamma_{0,s}\lambda^{-s-\frac{2}{3}}\,,
\end{align*}
as $\lambda\rightarrow\infty$ in $|\arg{}\lambda|\leq\frac{\pi}{2}-\delta<\frac{\pi}{2}$.  Recall that $\text{Ai}(z)$ is the Airy function (see Appendix \ref{sec:appendixc}).  In the above equations, $x$ is defined such that $\frac{4}{3}x^{3/2}=-2\zeta+3\log\left(\frac{2+e^{\zeta}}{2+e^{-\zeta}}\right)$, where $\zeta=\text{arccosh}\left(\frac{1}{4}z-1\right)$; this implies that $z>8\Leftrightarrow\zeta>0\Leftrightarrow{}x>0$. Additionally, the coefficients $\gamma_{j,s}$ are given (for $j=0,1$) by
\begin{align*}
	\gamma_{j,s}=\frac{(-1)^{s}s!}{(2\pi{}i)^{2}}\oint_{\{\text{sp}_{\pm}\}}\oint_{\{u(t)\}}\frac{t^{c-b-1}(t-1)^{a-c}(z+1-t)^{-a}u^{j}}{\left(\frac{1}{3}u^{3}-xu+\gamma-\ln\left(\frac{t-1}{t^{2}(z+1-t)}\right)\right)^{s+1}}{\rm d}u{\rm d}t\,,
\end{align*}
where $\gamma=-\frac{3}{2}\text{ln}(z+1)$. Here, the $t$-contour is a simple loop encircling the saddle points $\text{sp}_{\pm}:=2+e^{\pm\zeta}$, and, for each $t$ on this contour, the $u$-contour is a small loop encircling the point $u(t):=\left(\ln\left(\frac{t-1}{t^{2}(z+1-t)}\right)+xu-\gamma\right)^{1/3}$ \cite{OldeDaalhuis3}.
\end{itemize}

Using uniform asymptotic expansions to compute $\mathbf{F}$ is more helpful than doing so for $\mathbf{M}$ because the parameters are allowed to take complex values in the present expansions. The expansions for $\mathbf{F}$ can be used if $a$, $b$, $c$, and $\lambda$ are chosen so that the coefficients in the series expansions (which each entail either evaluating a contour integral or carrying out numerical differentiation) and the other special functions involved (which are detailed in Appendix \ref{sec:appendixc}) can be evaluated accurately. However, computing these functions is itself an expensive process, even if accurate results can be obtained, so we believe uniform asymptotic expansions are only likely to yield a potent numerical method if convergence is achieved after computation of very few terms of the series. We highlight that further research has been carried out on uniform asymptotic expansions for Gauss hypergeometric functions, for example in \cite{KOD,OldeDaalhuis4}.


\subsection{Other Methods for Computing $\mathbf{F}$} \label{sec:2F1_othermethods}

Another class of methods for computing $\mathbf{F}(a,b;c;z)$ or $_{2}F_{1}(a,b;c;z)$, which can be useful for $\left|z\right|<1$, is solving the \textbf{hypergeometric differential equation} \eqref{2f1de} numerically. To do this for $z$ close to the origin---provided that none of $c$, $c-a-b$, or $a-b$ are equal to an integer---we numerically integrate\eqref{2f1de} with initial conditions $f(0)=1$, $f'(0)=\frac{ab}{c}$, where the derivative at the origin follows from the Taylor series expansion \eqref{2F1}.   We use the code \url{d02bj} from the NAG Library \cite{NAGLibrary} for the numerical integration.  Integrating the differential equation is not limited by the radius of convergence of the Taylor series, hence this method is suitable for $|z| \geq 1$, provided that care is taken so that the integration contour avoids the branch point.  As for the confluent hypergeometric function, we highlight the relatively large computational cost of applying a differential equation based method, so such an approach should only be considered in parameter regimes where series methods fail.

Several other methods have been documented for the computation of $\mathbf{F}(a,b;c;z)$ \cite{Pearson}. For example, alternative methods can be used to compute the integral representation \eqref{2F1int}; these include splitting the integral, Romberg integration, and adaptive quadrature (though note that the latter two require that the integrand does not blow up at on either end point of the region of integration). Other possible methods include Pad\'{e} approximants \cite{PTVF}, evaluating a continued fraction representation \cite{Wolfram2} using the technique described in Ref.~\cite{AbramowitzStegun}, rational approximation \cite{Luke2}, Chebyshev expansion \cite{Luke3}, an $\epsilon$-expansion method \cite{HuangLiu}, exploiting alternative series expansions \cite{LopezTemme}, and developing relationships between different types of hypergeometric functions by evaluating Feynman path integrals \cite{KniehlTarasov}.


\subsection{Summary and Discussion} \label{sec:2F1_summary}

To compute $\mathbf{F}(a,b;c;z)$, we have considered several types of techniques---including Taylor series methods in Section \ref{sec:2F1_taylor}, the single-fraction method in Section \ref{sec:2F1_singlefraction}, and quadrature and differential equation methods in Sections \ref{sec:2F1_quadrature} and \ref{sec:2F1_othermethods}. In addition, we applied transformations and analytic continuation formulas from Sections \ref{sec:2F1_unitdisc}--\ref{sec:2F1_analyticcont} to find ways to compute $\mathbf{F}$ accurately and efficiently for all $z\in\mathbb{C}$.

We find that the series methods compute $\mathbf{F}(a,b;c;z)$ accurately for values of $\left|a\right|$ and $\left|b\right|$ less than about 50. The single-fraction method is particularly useful when $\left|c\right|<1$ and $\left|a\right|,\left|b\right|<30$. When ${\rm Re}(c)>{\rm Re}(b)>0$ or ${\rm Re}(c)>{\rm Re}(a)>0$, the method based on Gauss-Jacobi quadrature is often effective. A variety of the above methods work well if $\left|z\right|\lesssim 0.9$.

A difficulty arises when one needs to compute values of $\mathbf{F}(a,b;c;z)$ outside of the unit disk $\{|z| = 1\}$. In such situations, the transformation formulas of Section \ref{sec:2F1_unitdisc} or the methods detailed therein for the special cases satisfying $b-a\in\mathbb{Z}$ or $c-a-b\in\mathbb{Z}$ can be applied. A further issue arises when $\left|{\rm Re}(a)\right|$, $\left|{\rm Re}(b)\right|$, or $\left|{\rm Re}(c)\right|$ is too large for a method to work effectively on its own. (As a rough guide, this can occur when any of these values exceeds 50.) In this case, one can exploit the recurrence relation techniques of Section \ref{sec:2F1_largeparams_recurrences}.  One could alternatively employ the uniform asymptotic expansions of Section \ref{sec:2F1_largeparams_uaes} if $|a|$, $|b|$, or $|c|$ is large.

To devise an algorithm for cases in which all of the parameters and the variable are real, we follow the same type of procedure as for $\mathbf{M}$. Namely, we apply forward recurrences to $c$ so that all computations involve large values of $|c|$ and then use an additional transformation, such as \cite{DLMF}
\begin{align} \label{importanttrans2F1} 	&_{2}F_{1}(a,b;c;z)=(1-z)^{c-a-b}~_{2}F_{1}(c-a;c-b;c;z) \\
\ \nonumber \Leftrightarrow\quad&\mathbf{F}(a,b;c;z)=(1-z)^{c-a-b}~\mathbf{F}(c-a;c-b;c;z)\,.
\end{align}
We detail this strategy for the real case in Table \ref{2f1recommendationsreal}.

\begin{table}[h]
\renewcommand{\arraystretch}{1.15}
\begin{center}
\begin{footnotesize}
\begin{tabular}{|@{}c@{}|@{}c@{}|@{}l@{}|@{}c@{}|}
\hline
~Case~ & ~Regions for $a$, $b$, $z$~ & ~Recommended method(s) & ~Relevant sections and references~ \\ \hline \hline
I(A) & $a,b,c\geq0$, $z\geq0$ & ~Taylor series methods; & \ref{sec:2F1_taylor}, \cite{AbramowitzStegun,Muller} \\
 & & ~Single-fraction method; & \ref{sec:2F1_singlefraction} \cite{NPB1,NPB2} \\
 & & ~Gauss-Jacobi quadrature & \ref{sec:2F1_quadrature} \cite{Gautschi3} \\
 & & ~if ${\rm Re}(c)>{\rm Re}(b)>0$; & \\
 & & ~Adaptive ODE method; & \ref{sec:2F1_othermethods} \cite{NAGLibrary} \\ 
 & & ~Uniform asymptotic expansions & \ref{sec:2F1_largeparams_uaes} \cite{OldeDaalhuis2,OldeDaalhuis3,Temme7} \\
 & & ~for large parameter values & \\ \hline
I(B) & $a,b,c\geq0$, $z<0$ & ~Recurrences to combat cancellation, & \ref{sec:2F1_largeparams_recurrences} \cite{DeanoSegura,FLS,Gautschi1,GST1,Gautschi3,Temme1} \\
 & & ~then Taylor series & \ref{sec:2F1_taylor} \cite{AbramowitzStegun,Muller} \\
 & & ~or single-fraction method & \ref{sec:2F1_singlefraction} \cite{NPB1,NPB2} \\ \hline
II & $a,b<0$, $c\geq0$ & ~Transformation \eqref{importanttrans2F1}, & \ref{sec:2F1_summary} \cite{DLMF} \\
 & & ~then same as for Case I & \\ \hline
III & $a,b\geq0$, $c<0$ & ~Recurrence relations, & \ref{sec:2F1_largeparams_recurrences} \cite{DeanoSegura,FLS,Gautschi1,GST1,Gautschi3,Temme1} \\
 & & ~then Taylor series & \ref{sec:2F1_taylor} \cite{AbramowitzStegun,Muller} \\
 & & ~or single-fraction method & \ref{sec:2F1_singlefraction} \cite{NPB1,NPB2} \\ \hline
IV & $a,b,c<0$ & ~Transformation \eqref{importanttrans2F1}, & \ref{sec:2F1_summary} \cite{DLMF} \\
 & & ~then same as for Case III & \\ \hline
V & Either $a$ or $b<0$, & ~Recurrence relations, & \ref{sec:2F1_largeparams_recurrences} \cite{DeanoSegura,FLS,Gautschi1,GST1,Gautschi3,Temme1} \\
 & $c\geq0$ & ~then Taylor series & \ref{sec:2F1_taylor} \cite{AbramowitzStegun,Muller} \\
 & & ~or single-fraction method & \ref{sec:2F1_singlefraction} \cite{NPB1,NPB2} \\ \hline
VI & Either $a$ or $b<0$, & ~Recurrence relations, & \ref{sec:2F1_largeparams_recurrences} \cite{DeanoSegura,FLS,Gautschi1,GST1,Gautschi3,Temme1} \\
 & $c<0$ & ~then same as for Case V & \\ \hline
\end{tabular}
\caption{Recommendations for which methods to use for computation of the Gauss hypergeometric function for real parameters and a variable $z$ satisfying $|z|<1$. For $|z|\geq1$, appropriate transformations will also be necessary (see Section \ref{sec:2F1_unitdisc}).
}
\label{2f1recommendationsreal}
\end{footnotesize}
\end{center}
\end{table}

\begin{figure}
\centering
\includegraphics[bb=75 200 900 780,width=22cm]{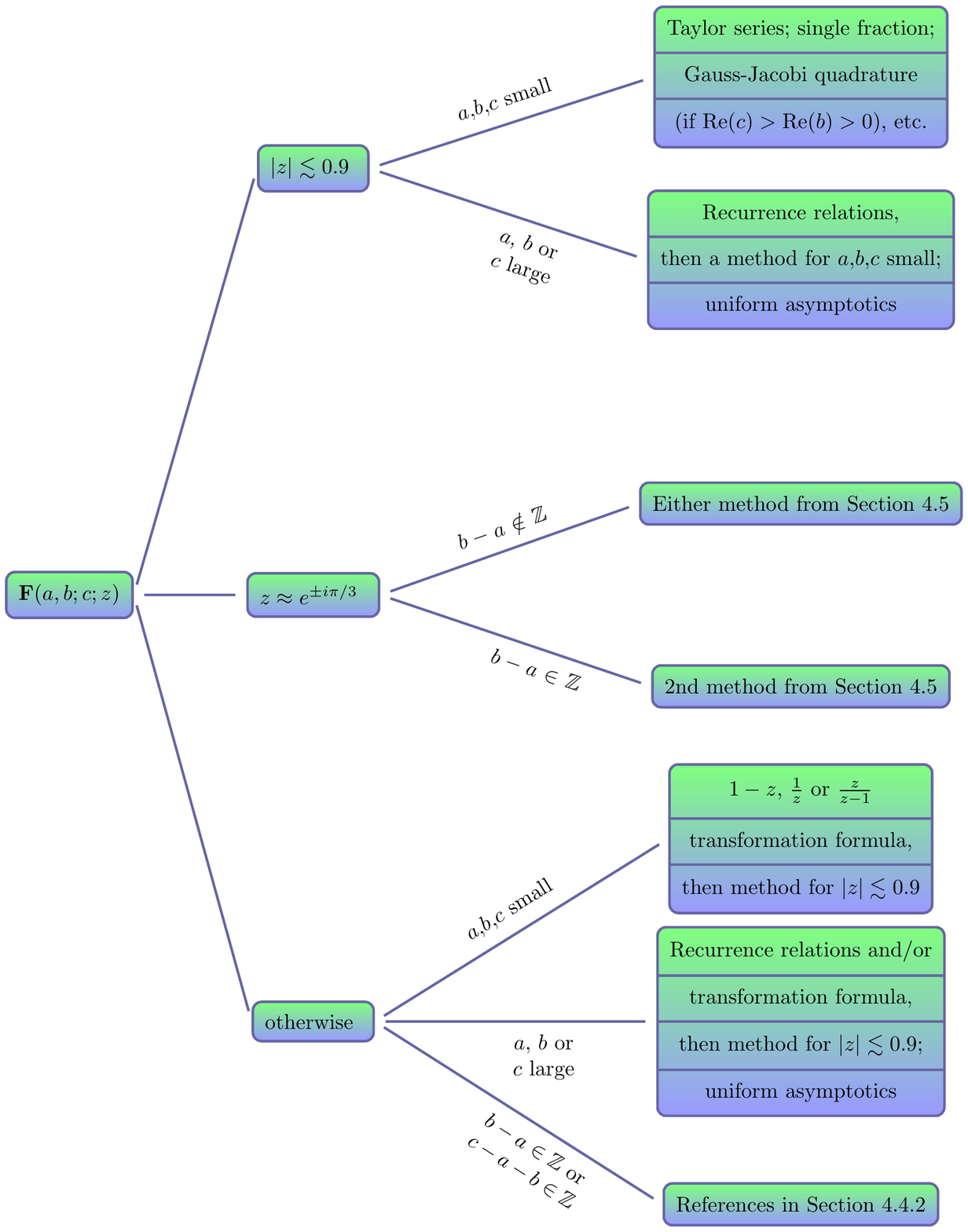}
\caption{Summary of recommended methods for computing $_{2}F_{1}$ for different values of the parameters $a$, $b$, and $c$ and the variable $z$.
}
\label{fig:summary2F1pic}
\end{figure}

When one extends the parameters and the variable into the complex plane, the problem of computing the Gauss hypergeometric function becomes more complicated.  One can achieve good results when computing $\mathbf{F}$ by following the procedure detailed in Fig.~\ref{fig:summary2F1pic}.

Details of other software written to compute the Gauss hypergeometric function can be found in Refs.~\cite{Forrey, Hsu, LozierOlver, MichelStoitsov, Moshier}.


\section{Concluding Remarks} \label{sec:conclusion}

The confluent and Gauss hypergeometric functions arise in a large and diverse set of applications, and many other functions are special cases of them.  However, both of these hypergeometric functions are notoriously difficult to compute effectively.

In this review paper, we have briefly outlined the theory of these hypergeometric functions, and we detailed a large number of good methods for computing them for different parameter and variable regimes. We have aimed to provide a roadmap for computing these functions for anyone who wishes to use them in a specific application and requires an effective implementation strategy for their evaluation.

There is no one method that provides a panacea for computing either $\mathbf{M}$ or $\mathbf{F}$, as there are numerous difficult computational issues (such as cancellation and overflow). Accordingly, we have detailed numerous good methods that should be chosen carefully depending on the values of the parameters $a$ and $b$ (and also $c$ for $\mathbf{F}$), and the variable $z$.  We have included roadmaps for the computation of both $\mathbf{M}$ and $\mathbf{F}$ that should be use with appropriate transformations and recurrence relations (which we have also detailed).  It is important to note that error bounds for the majority of these methods have not been widely researched (as most of the methods involve known series representations and hence can be calculated to arbitrary precision in infinite-precision arithmetic), but we believe the methods that we recommend have good error properties in the parameter and variable regions specified.

We have written and posted software for computing $\mathbf{M}(a;b;z)$ and $\mathbf{F}(a,b;c;z)$ \cite{PearsonCodeLink}. We find the methods that we include in the software for computing $\mathbf{M}$ to be effective for $|a|,|b|,|z|\lesssim100$ for real parameters and variables, and $|a|,|b|,|z|\lesssim70$ for complex variables and parameters; the methods for computing $\mathbf{F}$ are effective for $|a|,|b|,|c|\lesssim100$ and for any $z$ if the parameters and variables are real, and for $|a|,|b|,|c|\lesssim70$ and for any $z$ if they are complex. The potency of the methods is in spite of the fact that the code is implemented in double precision arithmetic, in contrast to the majority of widely-used software for hypergeometric functions. There are parameter and variable regimes in addition to these for which our software may be effectively applied---for instance for large $|z|$ when computing $\mathbf{M}$---and we are currently developing the software to expand the parameter and variable regimes for which it is 
effective.

\vspace{0.75em}


\textbf{Acknowledgements.}  

This research was conducted over a number of years, beginning as a project in the University of Oxford's MSc program in Mathematical Modelling and Scientific Computing so that the Numerical Algorithms Group (NAG) would be able to acquire enough information on the subject of computation of hypergeometric functions to achieve their goal of writing packages on the topic for the NAG Library \cite{NAGLibrary}.  The routines \texttt{s22ba} \& \texttt{s22bb} (for computing confluent hypergeometric functions), and  \texttt{s22be} \& \texttt{s22bf} (for Gauss hypergeometric functions), within the NAG Library, are substantially based on the work in this paper and the MSc project.  We thank the staff at NAG---especially Mick Pont, David Sayers, and Lawrence Mulholland---for productive discussions.  We also thank Andy Wathen, and Nick Trefethen for useful advice.  We also thank Frank Olver for providing us with an advanced copy of drafts of the chapters on hypergeometric functions in Ref.~\cite{DLMF}.  We are grateful to the many scientists who have been in contact about this work, and have brought to our attention new applications of hypergeometric functions.  JWP thanks the Engineering and Physical Sciences Research Council (EPSRC) and NAG for funding.


\appendix

\section{Test Cases for $\mathbf{M}$ and $\mathbf{F}$} \label{sec:appendixa}

This appendix details the test cases for $\mathbf{M}$ and $\mathbf{F}$ that we use in our numerical experiments in Appendix \ref{sec:appendixb}.

\begin{table}
\renewcommand{\arraystretch}{1.15}
\begin{center}
\begin{scriptsize}
\begin{tabular}{|c|c||c|c|}
\hline
Case & $(a,b,z)$ & Case & $(a,b,z)$ \\ \hline \hline
1 & $(0.1,0.2,0.5)$ & 21 & $(20,-10+10^{-9},-2.5)$ \\ \hline
2 & $(-0.1,0.2,0.5)$ & 22 & $(20,10-10^{-9},2.5)$ \\ \hline
3 & $(0.1,0.2,-0.5+i)$ & 23 & $(-20,-10+10^{-12},2.5)$ \\ \hline
4 & $(1+i,1+i,1-i)$ & 24 & $(50,10,200i)$ \\ \hline
5 & $(10^{-8},10^{-8},-10^{-10})$ & 25 & $(-5,(-5+10^{-9})(1+i),-1)$ \\ \hline
6 & $(10^{-8},10^{-12},-10^{-10}+10^{-12}i)$ & 26 & $(4,80,200)$ \\ \hline
7 & $(1,1,10+10^{-9}i)$ & 27 & $(-4,500,300)$ \\ \hline
8 & $(1,3,10)$ & 28 & $(5,0.1,-2+300i)$ \\ \hline
9 & $(500,511,10)$ & 29 & $(-5,0.1,2+300i)$ \\ \hline
10 & $(8.1,10.1,100)$ & 30 & $(2+8i,-150+i,150)$ \\ \hline
11 & $(1,2,600)$ & 31 & $(5,2,100-1000i)$ \\ \hline
12 & $(100,1.5,2.5)$ & 32 & $(-5,2,-100+1000i)$ \\ \hline
13 & $(-60,1,10)$ & 33 & $(-5,-2-i,1+(2-10^{-10})i)$ \\ \hline
14 & $(60,1,10)$ & 34 & $(1,10^{-12},1)$ \\ \hline
15 & $(60,1,-10)$ & 35 & $(10,10^{-12},10)$ \\ \hline
16 & $(-60,1,-10)$ & 36 & $(1,-1+10^{-12}i,1)$ \\ \hline
17 & $(1000,1,10^{-3})$ & 37 & $(1000,1,-1000)$ \\ \hline
18 & $(10^{-3},1,700)$ & 38 & $(-1000,1,1000)$ \\ \hline
19 & $(500,1,-5)$ & 39 & $(-10+500i,5i,10)$ \\ \hline
20 & $(-500,1,5)$ & 40 & $(20,10+1000i,-5)$ \\ \hline
\end{tabular}
\label{appatable1}
\caption{List of 40 test cases that we used to generate the results for $\mathbf{M}$ in Appendix \ref{sec:appendixb}. The justification for the selection of these cases is given in Ref.~\cite{Pearson}.}
\end{scriptsize}
\end{center}
\end{table}

\begin{table}
\renewcommand{\arraystretch}{1.15}
\begin{center}
\begin{scriptsize}
\begin{tabular}{|c|c||c|c|}
\hline
Case & $(a,b,c,z)$ & Case & $(a,b,c,z)$ \\ \hline \hline
1 & $(0.1,0.2,0.3,0.5)$ & 16 & $(-100,-200,-300+10^{-9},0.5\sqrt{2})$ \\ \hline
2 & $(-0.1,0.2,0.3,0.5)$ & 17 & $(300,10,5,0.5)$ \\ \hline
3 & $(0.1,0.2,-0.3,-0.5+0.5i)$ & 18 & $(5,-300,10,0.5)$ \\ \hline
4 & $(10^{-8},10^{-8},10^{-8},10^{-6})$ & 19 & $(10,5,-300.5,0.5)$ \\ \hline
5 & $(10^{-8},-10^{-6},10^{-12},-10^{-10}+10^{-12}i)$ & 20 & $(2+200i,5,10,0.6)$ \\ \hline
6 & $(1,10,1,0.5+10^{-9}i)$ & 21 & $(2+200i,5-100i,10+500i,0.8)$ \\ \hline
7 & $(1,-1+10^{-12}i,1,-0.8)$ & 22 & $(2,5,10-500i,-0.8)$ \\ \hline
8 & $(2+8i,3-5i,\sqrt{2}-\pi{}i,0.75)$ & 23 & $(2.25,3.75,-0.5,-1)$ \\ \hline
9 & $(100,200,350,i)$ & 24 & $(1,2,4+3i,0.6-0.8i)$ \\ \hline
10 & $(2+10^{-9},3,5,-0.75)$ & 25 & $(1,0.9,2,e^{i\pi/3})$ \\ \hline
11 & $(-2,-3,-5+10^{-9},0.5)$ & 26 & $(1,1,4,e^{i\pi/3})$ \\ \hline
12 & $(-1,-1.5,-2-10^{-15},0.5)$ & 27 & $(-1,0.9,2,e^{-i\pi/3})$ \\ \hline
13 & $(500,-500,500,0.75)$ & 28 & $(4,1.1,2,0.5+(0.5\sqrt{3}-0.01)i)$ \\ \hline
14 & $(500,500,500,0.75)$ & 29 & $(5,2.2,-2.5,0.49+0.5\sqrt{3}i)$ \\ \hline
15 & $(-1000,-2000,-4000.1,-0.5)$ & 30 & $(\frac{2}{3},1,\frac{4}{3},e^{i\pi/3})$ \\ \hline
\end{tabular}
\label{appatable2}
\caption{List of 30 test cases that we used to generate the results for $\mathbf{F}$ in Appendix \ref{sec:appendixb}. The justification for the selection of these cases is given in Ref.~\cite{Pearson}.}
\end{scriptsize}
\end{center}
\end{table}

In addition to using test cases such as these, other methods could also be used to test the accuracy of our computations. For instance, one could use tabulated values from sources such as Refs.~\cite{Slater1,Slater2,ZhangJin}, test against known relations to elementary or special functions \cite{AbramowitzStegun,BMOF,Luke4,DLMF}, or test against known recurrence relations \cite{AbramowitzStegun,DeanoTemme,GST3,DLMF} or Wronskians \cite{Luke3,Temme4,DLMF}.


\newpage 

\section{Table of Results for $\mathbf{M}$ and $\mathbf{F}$} \label{sec:appendixb}

\begin{table}
\renewcommand{\arraystretch}{1.15}
\begin{scriptsize}
\begin{center}
\begin{tabular}{|@{}c@{}||@{}c@{}|@{}c@{}|@{}c@{}|@{}c@{}|@{}c@{}|@{}c@{}|@{}c@{}|}
\hline
~Case~ & ~Taylor (a)~ & ~Taylor (b)~ & ~Single fraction~ & ~Buchholz~ & ~Asymptotic (a)~ & ~Asymptotic (b)~ & ~Gauss--Jacobi~ \\ \hline \hline
1 & 16 & 16 & 16 & A & A & A & 16 \\
2 & 16 & 16 & 16 & 16 & A & A & C \\
3 & 16 & 16 & 16 & A & A & A & 15 \\
4 & 16 & 16 & 15 & 15 & A & A & C \\
5 & 16 & 16 & 15 & 9 & A & A & C \\
6 & 8 & 15 & 15 & 0 & A & A & C \\
7 & 15 & 16 & 15 & 5 & A & A & C \\
8 & 15 & 15 & 15 & 15 & 16 & 16 & 14 \\
9 & 16 & 16 & 15 & A & A & A & B \\
10 & 15 & 16 & 0 & 13 & 15 & 15 & 13 \\
11 & A & A & 0 & A & 16 & 16 & 12 \\
12 & 16 & 15 & 15 & 15 & A & A & C \\
13 & 0 & 0 & 0 & 14 & A & A & C \\
14 & 15 & 15 & 15 & 15 & A & A & C \\
15 & 0 & 0 & 0 & 15 & A & A & C \\
16 & 16 & 16 & 16 & 15 & A & A & C \\
17 & 16 & 16 & 16 & 16 & A & A & C \\
18 & A & A & 0 & A & 16 & 16 & 13 \\
19 & 0 & 0 & 0 & 16 & A & A & C \\
20 & 0 & 0 & 0 & 15 & A & A & C \\
21 & 6 & 6 & 6 & 8 & A & A & C \\
22 & 16 & 16 & 15 & 8 & A & A & C \\
23 & 5 & 5 & 5 & 8 & A & A & C \\
24 & 0 & 0 & B & A & A & A & C \\
25 & 14 & 14 & 14 & 15 & A & A & C \\
26 & 15 & 16 & 0 & A & 15 & 15 & 0 \\
27 & 14 & 14 & 14 & A & A & A & C \\
28 & 0 & 0 & B & A & 14 & 14 & C \\
29 & 14 & 15 & 15 & A & 14 & 14 & C \\
30 & A & 0 & A & A & B & B & C \\
31 & A & A & B & A & 15 & 15 & C \\
32 & 15 & 15 & 15 & A & 15 & 15 & C \\
33 & 0 & 15 & 15 & A & A & A & C \\
34 & 4 & 15 & 16 & 12 & A & A & C \\
35 & 4 & 15 & 15 & 11 & A & A & C \\
36 & 15 & 15 & 15 & 0 & A & A & C \\
37 & A & A & 0 & A & A & A & C \\
38 & A & A & 0 & A & A & A & C \\
39 & A & 1 & A & 0 & A & A & C \\
40 & 16 & 16 & 15 & A & A & A & C \\
\hline
\end{tabular}
\caption{Table of results for various methods described in Section \ref{sec:1F1} for computing $\mathbf{M}$. We state the number of digits of accuracy obtained by applying each method to each test case from Appendix \ref{sec:appendixa}. We obtained these results using {\scshape Matlab} R2013a. We note again that label `A' indicates non-convergence of a series method, the label `B' indicates that the method was not effective because of overflow, and the label `C' indicates that the test case in question did not fall under the valid parameter regime for the Gauss-Jacobi quadrature method.
}
\label{appbtable1}
\end{center}
\end{scriptsize}
\end{table}

\begin{table}
\renewcommand{\arraystretch}{1.15}
\begin{scriptsize}
\begin{center}
\begin{tabular}{|@{}c@{}||@{}c@{}|@{}c@{}|@{}c@{}|@{}c@{}|@{}c@{}|}
\hline
~Case~ & ~Taylor (a)~ & ~Taylor (b)~ & ~Single fraction~ & ~Gauss--Jacobi~ & ~B\"{u}hring~ \\ \hline \hline
1 & 16 & 16 & 16 & 16 & A \\
2 & 16 & 16 & 16 & 15 & A \\
3 & 16 & 16 & 16 & C & 16 \\
4 & 15 & 15 & 15 & C & A \\
5 & 16 & 16 & 16 & C & 12 \\
6 & 15 & 15 & 15 & C & A \\
7 & 16 & 16 & 15 & C & A \\
8 & 14 & 15 & A & C & A \\
9 & A & A & A & B & A \\
10 & 16 & 16 & 15 & 16 & 0 \\
11 & 8 & 8 & 8 & C & A \\
12 & 16 & 16 & 16 & C & A \\
13 & 0 & 0 & 0 & C & A \\
14 & A & A & B & C & A \\
15 & 14 & 14 & 0 & C & A \\
16 & 0 & 0 & 0 & C & A \\
17 & A & A & 0 & C & A \\
18 & 0 & 0 & 0 & 13 & A \\
19 & 0 & 0 & 0 & C & A \\
20 & 0 & 0 & B & 6 & A \\
21 & 3 & 3 & A & B & A \\
22 & 16 & 16 & 16 & B & A \\
23 & A & A & 0 & C & 16 \\
24 & A & A & A & 0 & A \\
25 & A & A & 3 & 14 & 11 \\
26 & A & A & 5 & 14 & 11 \\
27 & 16 & 16 & 16 & 15 & 16 \\
28 & A & A & A & 14 & 16 \\
29 & A & A & A & C & 14 \\
30 & A & A & A & 15 & 16 \\
\hline
\end{tabular}
\caption{Table of results for various methods described in Section \ref{sec:2F1} for computing $\mathbf{F}$. We state the number of digits of accuracy obtained by applying each method to each test case stated in Appendix \ref{sec:appendixa}. We obtained these results using {\scshape Matlab} R2013a. We note again that label `A' indicates non-convergence of a series method, the label `B' indicates that the method was not effective because of overflow, and the label `C' indicates that the test case in question did not fall under the valid parameter regime for the Gauss-Jacobi quadrature method.
}
\label{appbtable2}
\end{center}
\end{scriptsize}
\end{table}

In this appendix, we present numerical results containing the number of digits of accuracy that we obtained when using a variety of methods for computing $\mathbf{M}$ and $\mathbf{F}$.

In Table \ref{appbtable1}, we show results that we obtained when computing $\mathbf{M}$ (or $_{1}F_{1}$ when $b$ is such that $\Gamma(b)$ is infinite in finite precision) for the two Taylor series methods of Section \ref{sec:1F1_taylor}, the single-fraction method of Section \ref{sec:1F1_singlefraction}, the Buchholz polynomial method of Section \ref{sec:1F1_buchholz}, the two methods for computing the asymptotic series of Section \ref{sec:1F1_asymptotic}, and the Gauss-Jacobi quadrature method of Section \ref{sec:1F1_quadrature}.

In Table \ref{appbtable2}, we present results that we obtained for computing $\mathbf{F}$ (or $_{2}F_{1}$ when $c$ is such that $\Gamma(c)$ is infinite in finite precision) using the two Taylor series methods of Section \ref{sec:1F1_taylor}, the single-fraction method of Section \ref{sec:2F1_singlefraction}, the Gauss-Jacobi quadrature method of Section \ref{sec:2F1_quadrature}, and the analytic continuation formula \eqref{continuation1} of Section \ref{sec:2F1_analyticcont}.

The label `A' in Tables \ref{appbtable1} and \ref{appbtable2} indicates that the series method in consideration had not converged after the computation of $500$ terms. The label `B' indicates that the method was not effective because of overflow. The label `C', which we use sometimes when we present a result from Gauss-Jacobi quadrature, indicates that the test case in question did not fall under the valid parameter regime for the method to be applicable.


\section{Other Special Functions Involved} \label{sec:appendixc}

To provide a comprehensive survey of the available methods for computing $\mathbf{M}$ and $\mathbf{F}$, we also need to consider the computation of other special functions that arise along the way. The most common such function is the Gamma function, which is required for several of the methods to compute $\mathbf{M}$ (see Sections \ref{sec:1F1_buchholz}, \ref{sec:1F1_asymptotic}, \ref{sec:1F1_quadrature}, \ref{sec:1F1_recurrences}, and \ref{sec:1F1_othermethods}) and $\mathbf{F}$ (see Sections \ref{sec:2F1_quadrature}, \ref{sec:2F1_unitdisc}, \ref{sec:2F1_analyticcont}, \ref{sec:2F1_largeparams_recurrences}, \ref{sec:2F1_largeparams_uaes}, and \ref{sec:2F1_othermethods}). Two other germane special functions are the incomplete gamma function, whose computation is necessary for the hyperasymptotic series in Section \ref{sec:1F1_asymptotic}, and the Bessel function, whose computation is necessary for the Buchholz polynomial method of Section \ref{sec:1F1_buchholz} and for the theory of uniform asymptotics. The 
psigamma,
Airy, and parabolic cylinder functions arise in the uniform asymptotic theory of Section \ref{sec:2F1_largeparams_uaes}.  These special functions are discussed in more detail in Refs.~\cite{GST3,Olver1}, and we note that the incomplete gamma and parabolic cylinder functions are themselves special cases of $\mathbf{M}$, the Airy function is a special case of $\mathbf{M}$, and the Bessel function is closely related to the \textbf{confluent hypergeometric limit function} $_{0}F_{1}$ \cite{DLMF}.\footnote{We define the confluent hypergeometric limit function $_{0}F_{1}$ as follows:
\begin{align*}
	_{0}F_{1}(-;a;z)=\sum_{j=0}^{\infty}\frac{1}{(a)_{j}}~\frac{z^{j}}{j!}\,.
\end{align*}}

The inbuilt {\scshape Matlab} routines for computing special functions (as of Version R2013a) are far from exhaustive; there only exist routines for the Gamma, Bessel, incomplete gamma, and Airy functions, and (apart from the Airy function) these are restrictive in terms of the the parameter and variable values that can be used. We therefore needed to use other routines to compute these functions accurately and efficiently for various parameter regions.

There is a inbuilt routine for computing the Gamma function $\Gamma(z)$ in {\scshape Matlab} (\url{gamma}), but it requires the argument $z$ to be real. The NAG Library \cite{NAGLibrary}, however, contains a routine (\url{s14ag}) to evaluate $\log[\Gamma(z)]$ for any complex $z$, and we used this for aspects of our work. If one wished to write his/her own routine for evaluating $\Gamma(z)$ for $z\in\mathbb{C}$, then for various values of $z$ one could use the Lanczos expansion \cite{Lanczos} (possibly using an expansion given by Godfrey at \cite{Godfrey}), Stirling expansions \cite{AbramowitzStegun}, Spouge's method \cite{Spouge}, or a Talbot contour method \cite{SchmelzerTrefethen,TWS,Weideman,WeidemanTrefethen}. When ${\rm Re}(z)<\frac{1}{2}$, the transformation \cite{PTVF}
\begin{align*}
	 \Gamma(z)\Gamma(1-z)=\frac{\pi}{\sin(\pi{}z)}\,,
\end{align*}
could be helpful when using the Talbot contour method.

{\scshape Matlab} also has an inbuilt routine to compute the Bessel function $J_{\nu}(z)$, but it is only applicable for real $\nu$, so it would be useful to design a routine that can also be used for complex $\nu$.  Such a routine could then be used to apply the methods in Sections \ref{sec:1F1_buchholz} and \ref{sec:2F1_largeparams_uaes}. In constructing such a code, one might wish to take advantage of its Taylor series, Hankel asymptotic representation, integral relations \cite{AbramowitzStegun}, as well as work on practical computation of this function in \cite{GST_Bessel}, and we exploit these methods in our code.

The `upper' and `lower' incomplete gamma functions $\Gamma(a,z)$ and $\gamma(a,z)$ can be computed for any complex $z$ using the inbuilt {\scshape Matlab} function \url{gammainc}, but this routine can only be applied for real $a\geq0$.  One would hence need to generalize this routine to implement the hyperasymptotic expansions detailed in Section \ref{sec:1F1_asymptotic}. To do this, one could exploit Taylor series methods \cite{PTVF}, asymptotic expansions \cite{AbramowitzStegun}, continued-fraction representations \cite{PTVF}, recurrence relations \cite{AbramowitzStegun}, and the relation \cite{PTVF}
\begin{align*}
	\gamma(a,z)+\Gamma(a,z)=\Gamma(a)\,.
\end{align*}
A computational procedure for the evaluation of $\Gamma(a,z)$ and $\gamma^{*}(a,z):=\frac{\gamma(a,z)}{\Gamma(a)}$ is discussed in Ref.~\cite{Gautschi2}.

The parabolic cylinder function $U(a,x)$ has yet to be implemented either in {\scshape Matlab} or in the NAG Library. A routine for computing this function could exploit Maclaurin series, asymptotic series, and recurrence relations \cite{GST3}; solving the underlying differential equation numerically \cite{AbramowitzStegun}; and other strategies \cite{Temme2}. We also highlight articles which consider the practical computation of parabolic cylinder functions in \cite{GST_PC1,GST_PC2}. For the Airy function $\text{Ai}(z)$ and its derivative $\text{Ai}'(z)$ (which are closely related to the Bessel function \cite{AbramowitzStegun}) and the psigamma function $\psi(z)=\frac{\Gamma'(z)}{\Gamma(z)}$, we use the NAG Library routines \url{s17dg} and \url{s14af}, respectively \cite{NAGLibrary}, and also highlight work carried out in \cite{GST_Airy}.


\section{Recurrence Relations Used in Section \ref{sec:2F1_largeparams_recurrences}} \label{sec:appendixd}

In this appendix, we state the recurrence relations for $_{2}F_{1}$ that we denoted in Section \ref{sec:2F1_largeparams_recurrences} as $(++0)$, $(00+)$, $(++-)$, and $(+0-)$.  In addition to the relations themselves, we also indicate their minimal solutions (see Table \ref{minsolns2f1}) that inspire these notations. We discussed their regions of validity in Section \ref{sec:2F1_largeparams_recurrences}.

\begin{footnotesize}
\begin{align*}
\ \underline{(++0):}~~&\textbf{Recurrence:}\quad(c-a-n)(c-b-n)(c-a-b-2n-1)y_{n-1} \\
\ &~~~~~~~~~~~~~~~~~~~~~~~~+(c-a-b-2n)\left[c(a+b-c+2n)+c-2(a+n)(b+n)\right. \\
\ &\left.~~~~~~~~~~~~~~~~~~~~~~~+z\{(a+b+2n)(c-a-b-2n)+2(a+n)(b+n)-c+1\}\right]y_{n} \\
&~~~~~~~~~~~~~~~~~~~~~~~~+(a+n)(b+n)(c-a-b-2n+1)(1-z)^{2}y_{n+1}=0\,, \\
\ \notag &\textbf{Solution:}\quad~~~~~y_{n}={}_{2}F_{1}(a+n;b+n;c;z)\,, \\
\ \underline{(00+):}~~&\textbf{Recurrence:}\quad(c+n)(c+n-1)(z-1)y_{n-1}+(c+n)\left[c+n-1\right. \\
\ &\left.~~~~~~~~~~~~~~~~~~~~~~~~-\{2(c+n)-a-b-1\}z\right]y_{n}+(c-a+n)(c-b+n)zy_{n+1}=0\,, \\
\ \notag &\textbf{Solution:}\quad~~~~~y_{n}={}_{2}F_{1}(a;b;c+n;z)\,, \\
\ \underline{(++-):}~~&\textbf{Recurrence:}\quad-(a-c+2n)(a-c+2n-1)(b-c+2n-1)(b-c+2n)zUy_{n-1} \\
&~~~~~~~~~~~~~~~~~~~~~~~~+(c-n)(c_{1}U+c_{2}V+c_{3}UV)y_{n} \\
\ &~~~~~~~~~~~~~~~~~~~~~~~~+(a+n)(b+n)(c-n)(c-n-1)(1-z)^{3}Vy_{n+1}=0\,, \\
\ &~~~~~~~~~~\text{where:}\quad c_{1}=(1-z)(b-c)(b-1)[a-1+z(b-c-1)]\,, \\
\ &~~~~~~~~~~~~~~~~~~~~~~c_{2}=b(b+1-c)(1-z)[a+z(b-c+2)]\,, \\
\ &~~~~~~~~~~~~~~~~~~~~~~c_{3}=c-2b-(a-b)z\,, \\
\ &~~~~~~~~~~~~~~~~~~~~~~U=z(a+b-c+1)(a+b-c+2)+ab(1-z)\,, \\
\ &~~~~~~~~~~~~~~~~~~~~~~V=(1-z)(1-a-b-ab)+z(a+b-c-1)(a+b-c-2)\,, \\
\ \notag &\textbf{Solution:}\quad~~~~~y_{n}={}_{2}F_{1}(a+n;b+n;c-n;z)\,, \\
\ \underline{(+0-):}~~&\textbf{Recurrence:}\quad{}z(a-c+2n)(a-c+2n-1)(b-c+n)[a+n+z(b-c+n+1)]y_{n-1} \\
\ &~~~~~~~~~~~~~~~~~~~~~~~~+(c-n)\left[(a+n-1)(c-n-1)(b-c+n)+(a+n)(a+n-1)\right. \\
\ &~~~~~~~~~~~~~~~~~~~~~~~~+(a+3b-4c+5n+2)z+(b-c+n)(b-c+n+1)(4a-c+5n-1)z^{2} \\
\ &~~~~~~~~~~~~~~~~~~~~~~~~-(a-b+n)(b-c+n)(b-c+n+1)z^{3}]y_{n} \\
\ &~~~~~~~~~~~~~~~~~~~~~~~~\left.-(a+n)(c-n)[a+n-1+z(b-c+n)\right](1-z)^{2}y_{n+1}=0\,, \\
\ \notag &\textbf{Solution:}\quad~~~~~y_{n}={}_{2}F_{1}(a+n;b;c-n;z)\,.
\end{align*}
\end{footnotesize}


\end{document}